\documentclass[aos]{imsart}
\RequirePackage{amsthm,amsmath,amsfonts,amssymb}
\RequirePackage[numbers]{natbib}
\usepackage{amssymb}

\startlocaldefs
\theoremstyle{plain}

\newtheorem{theorem}{Theorem}[section]
\newtheorem{lemma}[theorem]{Lemma}

\newtheorem{proposition}[theorem]{Proposition}
\newtheorem{conjecture}[theorem]{Conjecture}
\newtheorem{remark}[theorem]{Remark}

\newtheorem{corollary}[theorem]{Corollary}

\theoremstyle{remark}

\def \hat{\widehat}
\def \bar{\overline}
\newcommand{\eproof}{$\Box$}
\endlocaldefs

\begin{document}

\begin{frontmatter}
\title{Sampling and Filtering with Markov Chains
}
\runtitle{Sampling and Filtering with Markov Chains}

\begin{aug}
\author[A]{\fnms{Michael}~\snm{ A. Kouritzin
}\ead[label=e1]{michaelk@ualberta.ca}},
\address[A]{Department of Mathematical and Statistical Sciences,
University of Alberta\printead[presep={,\ }]{e1}}

\end{aug}

\begin{abstract}
A continuous-time Markov chain rate change formula for 
simulation, model selection, filtering and theory is proven.
It is used to develop Markov chain importance sampling, rejection
sampling, branching particle
filtering algorithms and filtering equations akin
to the Duncan-Mortensen-Zakai equation and the Fujisaki-Kallianpur-Kunita
equation but for Markov signals with general continuous-time Markov chain observations.
A direct method of solving these filtering equations is given that, for example, 
applies to trend, volatility and/or parameter estimation in
financial models given tick-by-tick
market data.
All the results also apply to continuous-time Hidden Markov Models (CTHMM), which have
become important in applications like disease progression tracking, as special cases
and the corresponding CTHMM results are stated as corollaries.
\end{abstract}



\end{frontmatter}

\maketitle

\setcounter{equation}{0}

\section{Introduction}\label{Intro}

Importance sampling (IS), introduced by Kloek and van Dijk \cite{KlvanD}, is an important
statistical variance reduction technique that is vital in problems
like Monte Carlo simulation (see e.g.\ \cite{KlvanDb}, \cite{ElMaLuBu}). 
From an alternative viewpoint,
IS is a general technique for estimating properties of a particular 
distribution, while only having samples generated from a different 
distribution.
Indeed, IS can also be important in rejection resampling and in
developing filtering equations, which will be further demonstrated
herein.
A method to extend IS to Markov chains (with possibly time and/or
hidden state dependent rates) will be introduced first in this paper.
Consider a process $Y$ (which will later be the observations) with finite or 
countable state space $O$ that is Markov.
Let it have a potentially time-dependent rate $\gamma_{i\rightarrow j}(s)$ of
going from $i\in O$ to $j\in O$ at time $s$.
Then, its generator is
\[
\mathbb L_sg(i)= \sum_{j\ne i,j\in O}\gamma_{i\rightarrow j}(s)
[g(j)-g(i)], 
\]
and $Y$ solves a well-posed martingale problem with $\mathbb L_s$ 
on a filtered probability space 
$(\Omega,\mathcal F, \{\mathcal F_t\}_{t\ge0}, P)$ such that
$\mathcal F^{Y}_{t}\doteq\sigma\{Y_s,s\le t\}\subset \mathcal F_t$ for all 
$t\ge 0$ under sufficient regularity. 
But, rather than simulating this chain directly, 
one samples $Y$ 
with respect to a reference measure $Q$ from
a simpler chain with generator: 
\begin{equation}\label{Lbar}
\overline {\mathbb{L}} g(i)= \sum_{j\ne i,j\in O}\overline \gamma_{i\rightarrow j}
[g(j)-g(i)],
\end{equation}
where the rates $\{\overline \gamma_{i\rightarrow j}\}_{i\ne j}$ may,
for example, not depend upon time or be more balanced.
Then, we will prove, under regularity conditions, that the likelihood ratio
process
\begin{equation}\label{FirstA}
A_{t}=\exp\left(\int _{0}^{t}\overline\gamma _{Y_{s\rightarrow} }-\gamma_{Y_{s\rightarrow} }\left({s} \right)ds\right)\prod\limits _{0<s\le t}\left[1+\left(\frac{\gamma_{Y_{s-}\rightarrow Y_{s}}\left({s} \right)}{\overline\gamma_{Y_{s-}\rightarrow Y_{s}}}-1\right)\Delta N_{s}\right]
\end{equation}
is a $\left\{\mathcal F^{Y}_{t}\right\}$-martingale under $Q$, where
$N$ counts the state changes of $Y$ and 
$\overline\gamma _{i\rightarrow }=\sum\limits_{j\ne i}\overline\gamma _{i\rightarrow j}$, $\gamma _{i\rightarrow }(s)=\sum\limits_{j\ne i}\gamma _{i\rightarrow j}(s)$. 
Moreover,
under the Girsanov-style measure change
$$\left.\frac{dP}{dQ}\right|_{\mathcal F_{t}}=A_{t},$$
$Y$ has the desired generator $\mathbb L_s$, hence desired process distribution, 
under $P$.  
This importance sampling result has important applications to 
things like simulation, model selection/verification, particle filters, 
parameter estimation, filtering theory and expanding the theory
of Hidden Markov models to Markov chain observations.
For example, we can simulate a Markov chain with (simple) proposal rates 
$\{\overline\gamma _{i\rightarrow j}\}$ and use rejection sampling
to create a Markov chain with (elaborate, possibly time-dependent) target rates 
$\{\gamma _{i\rightarrow j}(s)\}$.

Hidden Markov models (HMMs) were introduced in a series of papers by Baum and 
collaborators \cite{Baum1}, \cite{Baum2}.
Traditional HMMs have enjoyed termendous success in 
applications like computational finance \cite{Petro}, 
single-molecule kinetic 
analysis \cite{Nico}, speech recognition, and protein folding \cite{Stig}.
In such applications, the unobservable hidden states $X$ are a discrete-time or continuous-time Markov chain and the observations process $Y$ is some distorted, corrupted partial information or measurement of the current state of $X$ satisfying the condition
$$
P\left(Y_{t}\in A\big|X_{s},s\le t\right)=P\left(Y_{t}\in A\big|X_{t}\right)
$$
at update or observation times $t$.
These probabilities $P\left(Y_{t}\in A\big|X_{t}\right)$ are called the 
\emph{emission probabilities}.
In the case of a Continuous Time HMM (CTHMM), the times $\{t_i\}$ that the 
observations change are also random.
We consider the observations to be held constant in a cadlag manner until a new update occurs.
In this way, our results herein apply to CTHMM as corollaries.
CTHMM have been used successfully in applications like network
performance evaluation (see \cite{Weiwei}) and disease progression (see \cite{Liu2015}).

Still CTHMM's simple observations can be limiting.  
Consider observing weather conditions $Y$ based upon hidden climate state $X$.  
Wind, precipitation and temperature are not independent random samples but rather
time evolutions, depending upon the climate state.  
In other words, they depend upon their past values in addition to the current hidden state contradicting the Hidden Markov Model assumption.  
Perhaps, the simplest correction to the CTHMM would be to allow 
observation dependence only on their immediate past as well as the hidden state 
and insist that the measurements only take a finite or countable number of 
possibilities.
In this case, the observations become a Markov chain, whose transition probabilities depend upon the hidden states.  
However, for generality, we can allow the hidden state to be a general strong
Markov process on some state space $E$, by assuming:\vspace*{0.1cm}\\
({\bf C0}) $X$ is the unique
solution to a ($L,\mu$)-martingale problem (m.p.):
\[
M^{f}_{t}=f\left(X_{t} \right)-\int _{0}^{t}Lf\left(X_{s} \right)ds
\]
is a martingale for a linear operator $L$ and a rich enough class $D_L$ of bounded, continuous $f:E\rightarrow \mathbb R$ such that $\mathcal L\left(X_{0}\right)=\mu$.
(For convenience, we also let $\hat D_L=D_L\times \overline C (O)$.)\\
In the case that $X$ is actually a Markov chain, the Markov process 
infinitesimal generator $L$ for $x\in E$ is of the form
$Lf\left(x\right)=\sum\limits _{j\ne x}\lambda _{x\rightarrow j}
\left[f\left(j\right)-f\left(x\right)\right]$.
(The reader is referred to Ethier and Kurtz \cite{EtKu} for information
about the 
martingale problem and operator of a Markov process.
We will not require much of this theory other than noting $L1=0$
always for Markov generators and well-posednesss of martingale problems
is a standard way specifying models and defining process distributions.)

Given the hidden state, the observations can be modelled naturally as
a Markov chain depending upon this hidden state.
For example, given the hidden state of the climate, the observed 
wind, temperature etc.\
evolve as a Markov chain.
(In reality they would evolve as a Markov process, which could be fine
approximated by a Markov chain.)
In particular, we allow the observation rates to depend upon the hidden
state so the observations also satisfy a m.p. but with an 
operator that depends on the signal:
\begin{equation}\label{May21_22a}
m_{t}^{g}=g\left(Y_{t}\right)-\int _{0}^{t}\mathbb L\left(X_{s} \right)g\left(Y_{s}\right)ds
\end{equation}
is a martingale for all bounded, continuous $g$.
Here, 
$\mathbb L(x)g\left(y\right)=\sum\limits _{j\ne y}\gamma _{y\rightarrow j}(x)
\left[g\left(j\right)-g\left(y\right)\right]$, where 
$\{\gamma _{i\rightarrow j}(x)\}_{i,j\in O,i\ne j}$ are the rates
from state $i$ to state $j$ when the hidden state is $x$.
(This setting not only includes the nonlinear filtering setup but also 
things like common stochastic volatility models, both in the general 
Markov chain setting.)
If the observation noise is independent of the hidden state's, then
the role of independence gives us the combined m.p. for the
hidden state and observations:
\begin{equation}\label{Plus}
M^P_t=M^P_t(f,g)=f\left(X_{t} \right)g\left(Y_{t}\right)-\!\int _{0}^{t}g\left(Y_{s}\right)Lf\left(X_{s} \right)ds-\!\int _{0}^{t}f\left(X_{s} \right)\mathbb L\!\left(X_{s} \right)g\left(Y_{s}\right)ds\
\end{equation}
is a martingale for all bounded, continuous $g$ and $f\in D_L$.
For clarity, we will refer to the model of a hidden Markov signal $X$
(whether it is a chain or a general process) observed through a
continuous-time Markov chain $Y$ as a continuous Markov observation 
model (CMOM). 
It subsumes the popular CTHMM and can be studied through joint m.p. (\ref{Plus}).

A key observation for learning about the hidden state 
is that we can construct this model from a reference probability
$Q$ using importance sampling.
Motivated by the development of the Duncan-Mortensen-Zakai equation
(see \cite{Zakai}, \cite{KoLo08}), we can 
first consider simple, fake observations $Y$ with some reference probability $Q$
that do not depend upon the hidden state at all.
In particular, they could have rates $\{\overline \gamma_{i\rightarrow j}\}_{i\ne j}$
(that do not depend upon $x$) as above and be constructed to have the combined m.p.:
\begin{equation}\label{Star}
M^Q_t=f\left(X_{t} \right)g\left(Y_{t}\right)-\int _{0}^{t}g\left(Y_{s}\right)Lf\left(X_{s} \right)ds-\int _{0}^{t}f\left(X_{s} \right)\overline{\mathbb L}g\left(Y_{s}\right)ds
\end{equation}
is a martingale for all bounded, continuous $f,g$, where $\overline{\mathbb L}g\left(y\right)=\sum\limits _{j\ne y}\overline \gamma _{y\rightarrow j}\left[g\left(j\right)-g\left(y \right)\right]$. 
Then, adjusting likelihood ratio $A$ from (\ref{FirstA}) to account for the hidden state, one has that
\begin{equation}\label{AforX}
A_{t}=\exp\left(\int _{0}^{t}\overline\gamma _{Y_{s\rightarrow} }-\gamma_{Y_{s\rightarrow} }\left({X_s} \right)ds\right)\prod\limits _{0<s\le t}\left[1+\left(\frac{\gamma_{Y_{s-}\rightarrow Y_{s}}\left({X_s} \right)}{\overline\gamma_{Y_{s-}\rightarrow Y_{s}}}-1\right)\Delta N_{s}\right]
\end{equation}
is a $\left\{\mathcal F^{Y}_{t}\right\}$-martingale under $Q$ that
converts $Q$ into a new probability $P$, where the hidden state and observations $(X,Y)$ 
solve the desired joint m.p. (\ref{Plus}). 
(Here and below, $\gamma_{i\rightarrow}=\sum_{j\ne i}\gamma_{i\rightarrow j}$ is the rate
of leaving state $i$ and $N$ counts the state transitions of $Y$.)

The \emph{filter} $\pi _{t}\left(B\right)=P\left(X_{t}\in B\big|\mathcal F^{Y}_{t}\right)$, for Borel subsets $B$ of $E$, 
provides information on the hidden state based upon the model and the back
observations for both CTHMMs and CMOMs.  
$\{\pi_t,\ t\ge0\}$ is a probability measure-valued process.  
The \emph{unnormalized filter} $\sigma _{t}\left(B\right)=E^Q\left(A_t1_{X_{t}\in B}\big|\mathcal F^{Y}_{t}\right)$, 
with $A$ as defined in (\ref{AforX}), is a (finite, not-necessarily-probability) measure-valued process that
provides the filter through Bayes rule
\begin{equation}\label{Bayes}
\pi _{t}\left(f\right)=\frac{\sigma _{t}\left(f\right)}{\sigma _{t}\left(1\right)}
\end{equation}
for (at least) all bounded, measurable functions $f$, where
we used the notation $\pi _{t}\left(f\right)=\int _{E}fd\pi _{t}$
for an integrable function $f$.
However, the unnormalized filter provides more than a means to compute the filter.  
Rather as explained in \cite{KoZe05}, 
$\sigma _{t}\left(1\right)=\sigma _{t}\left(E\right)$ also provides 
the \emph{model rating} Bayes factor (i.e. integrated Likelihood)  of the model under consideration over the reference model.  
This is an overall rating of both overall structure and
model parameters, including those in the hidden state component of the model.   
If we have two sets of parameters, then we can form two models $M^1,M^2$ and produce two Bayes factors $\sigma^1_t(1),\sigma^2_t(1)$.  
We can compare these models on real data by evaluating $B_{1|2}(t)=\frac{\sigma ^{1}_{t}\left(1\right)}{\sigma ^{2}_{t}\left(1\right)}$.
These models' signals can be singular to each other, even of different
dimensions, so Bayes' factor methods are very general and effective.
One could even test if there is value in having a hidden state.
Hence, it is often more useful to produce a \emph{direct} or \emph{particle filter} approximation to the unnormalized filter than the (normalized)
filter.  

Particle filters that give Bayes factor information can be considered
\emph{model rating} particle filters.
Let $\left\{X^{i}\right\}_{i=1}^{\infty }$ be independent copies of the signal, 
called particles, recall the likelihood $A$ that was used to convert $Q$ into $P$, let
$$A^i_{t}=\exp\left(\int _{0}^{t}\overline\gamma _{Y_{s\rightarrow} }-\gamma_{Y_{s\rightarrow} }\left(X^i_{s} \right)ds\right)\prod\limits _{0<s\le t}\left[1+\left(\frac{\gamma_{Y_{s-}\rightarrow Y_{s}}\left(X^i_{s}\right)}{\overline\gamma_{Y_{s-}\rightarrow Y_{s}}}-1\right)\Delta N_{s}\right],
$$
which is like $A$ in (\ref{AforX}) except the hidden state is replaced with the
particle, and form
$$\sigma ^{N}_{t}\left(f\right)=\frac{1}{N}\sum\limits _{i=1}^{N}A^{i}_{t}f\left(X^{i}_t\right).$$
Then, by the independence of $Y$ and  $\{X^i\}_{i=1}^\infty$ under $Q$, we can fix the path $Y$ and find that 
\begin{equation}
\label{WeightPart}
\frac{1}{N}\sum\limits _{i=1}^{N}A^{i}_{t}f\left(X^{i}_t\right)\rightarrow E^Q\left[A_tf\left(X_{t}\right)\big|\mathcal F^{Y}_{t}\right]
\end{equation}
i.e. $\sigma ^{N}_{t}\left(f\right)\rightarrow \sigma_{t}\left(f\right)$ a.s. [$Q$] for each $f$ by the strong law of large numbers.  
By selecting a countable collection of $f$ that is closed under multiplication and strongly separate 
points, 
one can show a.s. convergence as measures i.e.  $\sigma ^{N}_{t}\left(\cdot\right)\Rightarrow \sigma_{t}\left(\cdot\right)$ a.s. [$Q$] 
(see Lemma 2 of \cite{BlKo10} and Lemma 7 of \cite{KoRe}). 
We refer to $\sigma^N_t$ as the weighted particle filter.  It approximates the unnormalized filter.
Notice that the real observations are used with each particle $X^i $ in $A^i$.

The problem with weighted particle filters is that the particles drift away from the signal over time (except when working on small compact sets) and do
not contribute much to the conditional distributional.
To maintain enough effective particles, we develop a branching method
analogous to that of \cite{Ko17a}.

The cellebrated Fujisaki-Kallianpur-Kunita (FKK) \cite{FKK72} and Duncan-Mortensen-Zakai
(DMZ) \cite{Zakai} equations were major breakthroughs in the classical nonlinear filtering problem yielding
the evolution of the filter and unnormalized filter.  
They are now the bases of many computational methods of solving filters
on real problems.
However, they are only known for the classical observation setting as well
as certain specific signal-dependent noise structures (see \cite{CrKoXi}) 
and not for general Markov chain observations.
As a final theoretical contribution we prove the FKK and DMZ equations for the
CMOM filtering problem, i.e. for Markov signals and continuous-time 
Markov chain observations.
In particular, for the setting described above, we show that
$\sigma $ is the
unique strong $D_{{\mathcal M}_f(E)}[0,\infty)$-valued solution to: 
\begin{eqnarray}\label{DMZsingf}
\!\!\sigma_t (f(\cdot))&\!=&\!\sigma_0 (f(\cdot))+\int_0^t\sigma_s (Lf(\cdot)) ds+\int_0^t\sigma_s (f(\cdot)(\overline\gamma_{Y_{s}\rightarrow }-{\gamma_{Y_{s}\rightarrow }\left(\cdot \right)})
ds\ \\\nonumber
&\!&\!+\int_{0}^t 
\sigma_{s-}\left(\left[f(\cdot )\frac{\gamma_{Y_{s-}\rightarrow Y_s}\left(\cdot  \right)}{\overline \gamma_{Y_{s-}\rightarrow Y_s}}-f(\cdot )\right]
\right)dN_s, \quad s.t.\ \sigma_0=\mathcal L(X_0).
\end{eqnarray}
for all $f\in D_L$ 
that also satisfies:
\begin{eqnarray}\label{DMZWea}
&&\!\!\sigma_t (f(\cdot,Y_t))\\\nonumber
&\!\!=&\sigma_0 (f(\cdot,Y_0))+\int_0^t\sigma_s (Lf(\cdot,Y_s)) ds+\int_0^t\sigma_s (f(\cdot,Y_{s})(\overline\gamma_{Y_{s}\rightarrow }-{\gamma_{Y_{s}\rightarrow }\left(\cdot \right)})
ds\ \ \ \\\nonumber
&\!&+\int_{0}^t 
\sigma_{s-}\left(\left[f(\cdot ,Y_s)\frac{\gamma_{Y_{s-}\rightarrow Y_s}\left(\cdot  \right)}{\overline \gamma_{Y_{s-}\rightarrow Y_s}}-f(\cdot ,Y_{s-})\right]
\right)dN_s
\end{eqnarray} 
for all $f\in\hat D_L$. 
Moreover, probability measure-valued process $\pi$ solves
\begin{eqnarray}\label{FKK}
\ \ \ \ \ \pi_t(f(\cdot,Y_t))&\!=&\!\pi_0 (f(\cdot,Y_0))
+\int_0^t\pi_s (Lf(\cdot,Y_{s})) ds \\\nonumber
&\!- &\!\int_0^t\!\pi_s(f(\cdot,Y_{s}) \gamma_{Y_{s}\rightarrow }\left(\cdot  \right)
)
-\pi_s (f(\cdot,Y_{s}))\pi_s (\gamma_{Y_{s}\rightarrow }\left(\cdot  \right)
)ds\\\nonumber
&\!+ &\!\!\int_0^t 
\frac {\pi_{s-} \left(
\gamma_{Y_{s-}\rightarrow Y_s}\left(\cdot  \right)f(\cdot,Y_s)\right)
-\pi_{s-} \!\left(\gamma_{Y_{s-}\rightarrow Y_s}\left(\cdot  \right) \right)\pi_{s-}(f(\cdot,Y_{s-}))}
{\pi_{s-}\left(\gamma_{Y_{s-}\rightarrow Y_s}\left(\cdot  \right)\right)}
dN_s,
\end{eqnarray}
for all $f\in\hat D_L$ subject to $\pi_0=\mathcal L(X_0)$.
In the CTHMM case, the observation rates do not depend the current
observation.
Instead, we have a rate for updates $\gamma(x)$ that can depend upon
the hidden state and an emission probability $q_{y'}(x)=P(Y_t=y'|X_t=x,W_n=t)$
for all $t>0$, where $W_n$ is the time of the $n^{th}$ emission so 
$\overline \gamma_{i\rightarrow j}$, $\gamma_{i\rightarrow j}(x)$ become $\overline \gamma\, \overline q_{ j}$, $\gamma(x)q_{ j}(x)$,
where $\overline \gamma$ and $\overline q_{j}$ are some canonical
update rate and emission probability mass function that do not depend upon the hidden state
$x$.
In this CTHMM case, (\ref{Star}), (\ref{Plus}) and (\ref{DMZsingf}) become
\begin{equation}\label{CTHMMStar}
\!\!M^Q_t=f\left(X_{t} ,Y_{t}\right)-\int _{0}^{t}Lf\left(X_{s},Y_{s} \right)ds-\int _{0}^{t}\sum_j\overline\gamma\,
\overline q_{ j}[f\left(X_{s},j \right)-f\left(X_{s},Y_{s}\right)]ds,
\end{equation}
\begin{equation}\label{CTHMMPlus}
M^P_t=f\left(X_{t} ,Y_{t}\right)-\!\int _{0}^{t}Lf\left(X_{s},Y_{s} \right)ds-\!\int _{0}^{t}\sum_j\gamma(X_{s} )\,
 q_{j\,}\!(X_{s} )[f(X_{s} ,j)-f(X_{s} ,Y_{s})]ds
\end{equation}
for all $f\in\hat D_L$, where $L$ operates only on the first variable, and
\begin{eqnarray}\label{CTHMMDMZ}
\!\!\sigma_t (f(\cdot))&\!=&\!\sigma_0 (f(\cdot))+\int_0^t\sigma_s (Lf(\cdot)) ds+\int_0^t\sigma_s (f(\cdot)(\overline\gamma -{\gamma\left(\cdot \right)})
ds\ \\\nonumber
&\!&\!+\int_{0}^t 
\sigma_{s-}\left(\left[f(\cdot )\frac{\gamma(\cdot )q_{ Y_s}(\cdot )}{\overline \gamma\, \overline q_{Y_s}}-f(\cdot )\right]
\right)dN_s, \quad s.t.\ \sigma_0=\mathcal L(X_0).
\end{eqnarray}
for all $f\in D_L$.
(\ref{DMZWea}, \ref{FKK}) also simplify accordingly (see below).
(\ref{CTHMMDMZ}) can be used for filtering, model selection and as
a basis for parameter estimation in CTHMM.

Finally, we develop a direct solution approach solving 
(\ref{DMZsingf}) (or (\ref{CTHMMDMZ})) that is being used in current work to estimate trends and
volatility in financial models based upon tick-by-tick data.
In another application, the Bayes factor developed herein
is being used to identify the parameters and compare fitness of
the standard stochastic volatility models (with Markov chain
approximation applied to the price in order to match actual prices) to
tick-by-tick data.
Naturally, there are many other potential applications
for our results. 

\subsection{Layout}
The next section is a proof of our rate-change Girsanov's theorem.
Section 3 has our applications to resampling Markov chains, 
Monte Carlo simulation and particle
filtering and model selection with Markov chain observations.
Section 4 contains a new development for
filtering problems with Markov chain observations, establishing 
FKK and DMZ type equations for these filtering problems.
It also contains our direct  
solutions to these filtering equations.
Finally, our conclusions and highlights are in Section 5.
Due to space consideration, our applications are left to future work. 

\subsection{Notation}
Let
$B(S)$, $C(S)$ and $\bar{C}(S)$ be the bounded, continuous and
continuous bounded real functionals respectively  
and $M(S,S)$ be the measurable functions from $S$ to
$S$ with the topology of pointwise convergence
on any Polish space $S$.
Further, we let ${\mathcal M}_f(S)$ be the space of
finite (non-negative) Borel measures with the topology of weak convergence 
i.e. the notion that $\nu^n \Rightarrow\nu$
if and only if $\nu^n(f)\rightarrow\nu(f)$ for all $f\in\bar{C}(S)$, and
${\mathcal P}(S)\subset {\mathcal M}_f(S)$ be the probability measures. 
Finally, we let $D_S[0,\infty)$ denote the space cadlag path from
$[0,\infty)$ to $S$ equipped with the Skorokhod metric.

For $z>0$, $\lfloor z\rfloor$ is the greatest integer not more than
$z$ and $\{z\}=z-\lfloor z\rfloor$.

$\lambda_{i\rightarrow j}$ denotes the rate of a Markov chain going from state
$i$ to state $j$ and $\lambda_{i\rightarrow }=\sum_{j\ne i}\lambda_{i\rightarrow j}$ denotes the rate of leaving state $i$.

$a_{i,k}\stackrel{i}\ll b_{i,k}$ means $\forall k$, $\exists c_k > 0$ not depending on i s.t. $|a_{i,k}|\le c_k |b_{i,k}|$ $\forall i, k$.

$L$ will be used for a generator, a linear operator $C\left(S\right)\rightarrow C\left(S\right)$.
Here, $S$ will either be the hidden state space $E$ or the observation state space
$O$.

$A,A^i$ will be used for likelihood, particle weight processes respectively.

$\mathcal L\left(Z \right)$ will be used to denote the law or distribution of random variable $Z$.

$\delta_x$ will be used for Dirac delta measure at point $x$.

\section{Rate Change Formula}

In this section, we derive our rate-change Girsanov type theorem
for Markov chains and hidden Markov signals subject to Markov chain observations.
We proceed in the notation of the later application as the reader can 
just take $X_s=s$, $E=[0,T]$ for the former.

The following conditions will be imposed in our first main result to follow:
\begin{itemize}
\item[(C1)] The observation state space $O$ is a finite or countable space and $\sup\limits_{i\in O}\overline\gamma_{i\rightarrow }<\infty$.
\item[(C2)] $\sup\limits_{x\in E,i\in O}\frac{\gamma_{i\rightarrow }(x)}{\overline\gamma_{i\rightarrow }}<\infty$.
\item[(C3)] There are no cemetery states, meaning $\gamma_{i\rightarrow }(x), \overline\gamma_{i\rightarrow }>0$ for all $i\in O$, $x\in E$.
\end{itemize}

Under (C1), $\overline {\mathbb L}$, defined in (\ref{Lbar}), is a bounded operator and, given a distribution for $Y_0$ on $O$, there is a unique
solution to the m.p.:
\[
g(Y_t)-\int_0^t \overline {\mathbb L} g(Y_s)ds
\]
is a martingale for all continuous, bounded $g$ which can be constructed and simulated
in the following manner.
On our reference measure $Q$,
let $\left\{T^{i\rightarrow }_{n}\right\}_{n\in \mathbb{N},i\in O}$ be independent exponential random variables with rates  $\left\{\overline \gamma_{i\rightarrow }\right\}_{i\in O}$, 
let $\left\{\xi ^{i}_{n},n\in \mathbb{N}_0\right\}_{i\in O}$ be independent 
discrete random variables such that
$P\left(\xi ^{i}_{n}=j\right)=\frac{\overline{\lambda}_{i\rightarrow j}}{\overline{\lambda}_{i\rightarrow }}$ for $j\ne i\in O$ and let $\theta_0$ be a random sample from
$\mathcal L(Y_0)$.
Take all these to be independent of each other and of the hidden state $X$.
Then, $Y$ could be constructed under $Q$ in three steps starting from some $Y_{0}=\theta _{0}$, $W_0=0$: 

1) $\theta _{n}=\xi ^{\theta _{n-1}}_{n-1}$ for $n \in \mathbb N$ (gets the transitions as a discrete chain) 

2) $W_{n}=W_{n-1}+T_{n}^{\theta _{n-1}\rightarrow }$ for $n \in \mathbb N$ (get the transition times)

3) $Y_{t}=\theta _{n}$ for $t\in \left[W_{n},W_{n+1}\right)$ $n \in \mathbb N_0$ (create continuous time chain). 

Notice that the quadratic variation $\left[f\left(X\right),g\left(Y\right)\right]=0$ for all $f,g$ by 
the $X,Y$ independence under $Q$ so $\left[f\left(X\right),N\right]=0$,
where $N$ counts the transitions in $Y$, for all 
$f$.

Now, we have our first main result, which is a Girsanov change-of-measure
result changing Markov chains by
weighting.
\begin{theorem}\label{ImportanceTheorem}
Suppose (C1, C2, C3) hold and $\left(X,Y\right)$ satisfies the (\ref{Star}) martingale problem starting from
some initial law $\mathcal L(X_0,Y_0)=\nu$ under $Q$.
Then, 
\begin{equation}\label{Aweight}
A_{t}=\exp\left(\int _{0}^{t}\overline\gamma _{Y_{s}\rightarrow }-\gamma_{Y_{s} \rightarrow}\left(X_{s} \right)ds\right)\prod\limits _{0<s\le t}\left[1+\left(\frac{\gamma_{Y_{s-}\rightarrow Y_{s}}\left(X_{s} \right)}{\overline\gamma_{Y_{s-}\rightarrow Y_{s}}}-1\right)\Delta N_{s}\right]
\end{equation}
is a $\left\{\mathcal F^{Y}_{t}\right\}$-martingale under $Q$, where $N$ counts
the transitions of $Y$. 
Moreover, if we define a new probability measure via
$
\left.\frac{dP}{dQ}\right|_{\mathcal F_{t}}=A_{t},\ \ \forall t\ge 0,
$
then $(X,Y)$ satisfies the (\ref{Plus}) m.p. starting from
$\mathcal L(X_0,Y_0)=\nu$ under $P$.  
\end{theorem}

{\bf Note:} $\prod\limits _{0<s\le t}\left[1+\left(\frac{\gamma_{Y_{s-}\rightarrow Y_{s}}\left(X_{s} \right)}{\overline\gamma_{Y_{s-}\rightarrow Y_{s}}}-1\right)\Delta N_{s}\right]$ is $1$ except at the transition times $s$ at which times it is $\frac{\gamma_{Y_{s-}\rightarrow Y_{s}}\left(X_{s}\right)}{\overline\gamma_{Y_{s-}\rightarrow Y_{s}}}$.

\proof
By our construction, $N$ counts the transitions of $Y$ and the
$\{\xi_i^k\}$ determine the actual transitions.
Let $\mathcal F_{t}\doteq\sigma\{N_s,\xi_i^k,X_u:\,i\le N_s ,k\in O,s\le t,u\in[0,\infty)\}$, 
$\overline q_{i\rightarrow j}=\frac{\overline \lambda_{i\rightarrow j}}{\overline \lambda_{i\rightarrow }}$ and
$q_{i\rightarrow j}(x)=\frac{\lambda_{i\rightarrow j}(x)}{\lambda_{i\rightarrow }(x)}$.
Under $Q$, the combined chain $\left(Y,N\right)$ satisfies the m.p.
\begin{equation}\label{May30_22c}
g\left(Y_{t},N_{t}\right)-\int _{0}^{t}\overline {\mathbb L}^{N}g\left(Y_{s},N_{s}\right)ds
\end{equation}
is a $\{\mathcal F_{t}\}$-martingale, 
where  $\overline{\mathbb L}^{N}g\left(y,n\right)=\sum\limits _{j\ne y}\overline \gamma _{y\rightarrow j}\left[g\left(j,n+1\right)-g\left(y ,n\right)\right]$.
In particular, taking $g\left(y,n\right)=n$
\begin{equation}\label{MNdef}
M^N_t=N_t-\int _{0}^{t}\overline\gamma _{Y_{s\rightarrow}}ds
\end{equation}
and so
\begin{eqnarray}\label{exclaim1}
&&\int_0^t\! A_{s-}^{1}\!\left(\frac{\gamma_{Y_{s-}\rightarrow }\left(X_{s} \right)}{\overline\gamma_{Y_{s-}\rightarrow }}-1\!\right)\!dM^N_s\\\nonumber
&\!\!=&\!\!\int_0^t\! A_{s-}^{1}\!\left(\frac{\gamma_{Y_{s-}\rightarrow }\left(X_{s} \right)}{\overline\gamma_{Y_{s-}\rightarrow }}-1\!\right)\!dN_s-\int _{0}^{t}
\!A_{s-}^{1}\!\left(\frac{\gamma_{Y_{s-}\rightarrow }\left(X_{s} \right)}{\overline\gamma_{Y_{s-}\rightarrow }}-1\!\right)
\overline\gamma _{Y_{s\rightarrow}}ds,
\end{eqnarray}
with $A_{t}^{1}$ defined below, are $Q$-local martingale.
We break this weight $A$ from (\ref{Aweight}) into two factors $A^1$ and $A^2$ defined by $A^1_0=A^2_0=1$ and
\begin{eqnarray}\label{A1diffdef}
dA^{1}_{t}&=&A_{t-}^{1}\left[\left(\frac{\gamma_{Y_{t-}\rightarrow }\left(X_{t} \right)}{\overline\gamma_{Y_{t-}\rightarrow }}-1\right)d N_{t}+[\overline\gamma _{Y_{t\rightarrow} }-\gamma_{Y_{t\rightarrow} }\left(X_{t} \right)]dt\right]\\\label{A2diffdef}
dA^2_t&=&A_{t-}^{2}\left(\frac{q_{Y_{t-}\rightarrow Y_t}\left(X_{t} \right)}{\overline q_{Y_{t-}\rightarrow Y_t }}-1\right)d N_{t}.
\end{eqnarray}
To show this, we use integration by parts with $\alpha=A^1 A^2$
\begin{eqnarray}
&&d\alpha _{t}\label{May10_22b}=
\alpha _{t-}\left[\left(\frac{\gamma_{Y_{t-}\rightarrow }\left(X_{t} \right)}{\overline\gamma_{Y_{t-}\rightarrow }}-1\right)d N_{t}+[\overline\gamma _{Y_{t\rightarrow} }-\gamma_{Y_{t\rightarrow} }\left(X_{t}\right)]dt\right]\\
&+&\nonumber 
\alpha _{t-}\left(\frac{q_{Y_{t-}\rightarrow Y_t}\left(X_{t} \right)}{\overline q_{Y_{t-}\rightarrow Y_t }}-1\right)d N_{t}+ \alpha _{t-}\left(\frac{\gamma_{Y_{t-}\rightarrow }\left(X_{t} \right)}{\overline\gamma_{Y_{t-}\rightarrow }}-1\right)\left(\frac{q_{Y_{t-}\rightarrow Y_t}\left(X_{t}\right)}{\overline q_{Y_{t-}\rightarrow Y_t }}-1\right)d N_{t}\\
&=& \nonumber\alpha _{t-}\left(\frac{\gamma_{Y_{t-}\rightarrow }\left(X_{t}\right)}{\overline\gamma_{Y_{t-}\rightarrow }}\frac{q_{Y_{t-}\rightarrow Y_t}\left(X_{t} \right)}{\overline q_{Y_{t-}\rightarrow Y_t }}-1\right)d N_{t}+ \alpha _{t-}[\overline\gamma _{Y_{t\rightarrow} }-\gamma_{Y_{t\rightarrow} }\left(X_{t} \right)]dt\\
&=&\nonumber \alpha _{t-}\left(\frac{\gamma_{Y_{t-}\rightarrow Y_t}\left(X_{t} \right)}{\overline \gamma_{Y_{t-}\rightarrow Y_t }}-1\right)d N_{t}+ \alpha _{t-}[\overline\gamma _{Y_{t\rightarrow} }-\gamma_{Y_{t\rightarrow} }\left(X_{t} \right)]dt.
\end{eqnarray}
But, this stochastic exponential equation has unique solution (see Protter \cite[Theorem II.36]{P90})
$$\alpha_t=\exp\left(\int _{0}^{t}\overline\gamma _{Y_{s\rightarrow} }-\gamma_{Y_{s\rightarrow} }\left(X_{s} \right)ds\right)\prod\limits _{0<s\le t}\left[1+\left(\frac{\gamma_{Y_{s-}\rightarrow Y_{s}}\left(X_{s} \right)}{\overline\gamma_{Y_{s-}\rightarrow Y_{s}}}-1\right)\Delta N_{s}\right]=A_{t}$$
so $A_t=A^1 _t A^2_t$.  
Moreover, by Ito's formula (see Theorem II.36 of Protter) again
\begin{equation}\label{May20_22b}
A^{2}_{t}=\prod\limits _{0<s\le t}\left[1+\left(\frac{q_{Y_{s-}\rightarrow Y_{s}}\left(X_{s} \right)}{\overline q_{Y_{s-}\rightarrow Y_{s}}}-1\right)\Delta N_{s}\right]
\end{equation}
and
\begin{eqnarray}\label{A1def}
A^1_t&=&\exp \left(\int _{0}^{t}\overline\gamma _{Y_{s\rightarrow} }
-\gamma_{Y_{s\rightarrow} }\left(X_{s} \right)ds\right)\prod\limits _{0<s\le t}\left[1+\left(\frac{\gamma_{Y_{s-}\rightarrow }\left(X_{s}\right)}
{\overline\gamma_{Y_{s-}\rightarrow }}-1\right)\Delta N_{s}\right]\\
&=&\exp \left(\int _{0}^{t}\left[\overline\gamma _{Y_{s\rightarrow} }-\gamma_{Y_{s\rightarrow} }\left(X_{s} \right)ds+\ln\left(\frac{\gamma_{Y_{s-}\rightarrow }\left(X_{s} \right)}{\overline\gamma_{Y_{s-}\rightarrow }}\right)dN_{s}\right]\right).\nonumber
\end{eqnarray}
(\emph{Show} $A^1$ \emph{is a martingale}) 
By (\ref{exclaim1},\ref{A1diffdef}) $A^1$ is a $Q$ local martingale.
Under (C1, C2), it is a martingale by (\ref{A1def}).
Let $P^1$ be defined by 
$$\left.\frac{dP^1}{dQ}\right|_{\mathcal F_{t}}=A^1_{t}, \ \forall t\ge 0.$$
(\emph{Apply Girsanov-Meyer}) 
It follows by (\ref{Star}) with $f\equiv 1$ and Theorem III.20 of Protter that
\begin{equation}\label{May20_22c}
\widetilde{M}_t=g(Y_t)-g(Y_0)-\int _{0}^{t}\overline {\mathbb L}g(Y_s)ds-\int_0^t\frac1{A^1_s}d[A^1,M^Q]_s
\end{equation}
is a local martingale under $P^1$. 
However, by (\ref{Star}) with $f\equiv 1$ again, (\ref{A1diffdef}), (\ref{A1def})
\begin{eqnarray}
\left[A^{1},M^{Q}\right]_{t}&\!=&\!\left[A^{1},g(Y)\right]_{t}\\\label{May20_22d}
\ \ \ \int_0^t\!\frac1{A^1_s}d[A^1,M^Q]_s&\!=&\!\int
_{0}^{t}\frac{A^1_{s-}}{A^1_s}\left(\frac{\gamma_{Y_{s-}\rightarrow}\left(X_{s}
\right)}{\overline\gamma_{Y_{s-}\rightarrow }}-1\right)\left[g\left(\xi
^{Y_{s-}}_{N_{s-}}\right)-g\left(Y_{s-}\right)\right]d N_{s}\\\nonumber
&\!=&\!\int _{0}^{t}\!\frac{\overline\gamma_{Y_{s-}\rightarrow }}{\gamma_{Y_{s-}\rightarrow }\left(X_{s}\right)}\left(\frac{\gamma_{Y_{s-}\rightarrow }\left(X_{s} \right)}{\overline\gamma_{Y_{s-}\rightarrow }}-1\right)\!\left[g\!\left(\!\xi ^{Y_{s-}}_{N_{s-}}\right)-g\left(Y_{s-}\right)\right]\!d N_{s}.
\end{eqnarray}
Putting (\ref{May20_22c}, \ref{May20_22d}) together, we have that
$$\widetilde M_t=\int _{0}^{t}\frac{\overline\gamma_{Y_{s-}\rightarrow }}{\gamma_{Y_{s-}\rightarrow }\left(X_{s}\right)}\left[g\left(\xi ^{Y_{s-}}_{N_{s-}}\right)-g\left(Y_{s-}\right)\right]d N_{s}-\int _{0}^{t}\overline {\mathbb L} g(Y_{s})ds$$
and so 
$$\int _{0}^{t}\frac{\gamma_{Y_{s-}\rightarrow }\left(X_{s} \right)}{\overline\gamma_{Y_{s-}\rightarrow }}d\widetilde M_{s}=g\left({Y_{t}}\right)-g\left(Y_0\right)-\int _{0}^{t}\frac{\gamma_{Y_{s-}\rightarrow }\left(X_{s} \right)}{\overline\gamma_{Y_{s-}\rightarrow }}\overline {\mathbb L} g(Y_{s})ds$$
is a $P^1$-local martingale.\\
(\emph{New $Y$ mp}) Rewriting this and using (C1-C3) as well as bounded $g$, we have that
\begin{equation}\label{May30_22b}
g\left({Y_{t}}\right)-g\left(Y_0\right)-\int _{0}^{t}{\gamma_{Y_{s}\rightarrow }\left(X_{s} \right)}\sum_{j\ne Y_s}{\overline q_{Y_{s}\rightarrow j}}[g(j)- g(Y_{s})]ds
\end{equation}
is a $P^1$ martingale and following (\ref{May30_22c},\ref{MNdef}) that
\begin{equation}
N_{t}-\int _{0}^{t}\gamma_{Y_{s-}\rightarrow }\left(X_{s} \right)ds
\label{May30_22d}
\end{equation}
is a $P^1$-local martingale.\\
(\emph{Show} $A^2$ \emph{is a $P^1$-martingale}) 
Under (C1,C2), $\mathfrak L(x)g(y)\doteq \gamma_{y\rightarrow }(x)\sum_{j\neq y}
\overline q_{y\rightarrow j}[g(j)-g(y)]$ are uniformly (in $x$) bounded operators
and for a given $\{X_t,\,t\ge 0\}$, the $\mathfrak L(X_t)$ m.p.
is well posed.
The solution can be constructed using the same 3 step procedure given above except
in 2) the distribution of $T_n$ is determined by 
$P(T_n^{\theta_{n-1}}>t)=
e^{-\int_{W_n}^{W_n+t}\gamma_{\theta_{n-1}\rightarrow}(X_s)ds}$ for $t\ge 0$.
Hence, we still use independent $\{\xi^k_n\}$ such that $P^1\left(\xi ^{k}_n=j\right)=Q\left(\xi ^{k}_{n}=j\right)=\overline q_{k\rightarrow j}$
and
$$Y_s=\left\{\begin{array}{cc}\xi^{Y_{s-}}_{N_{s-}}&\Delta N_s=1\\Y_{s-}&\Delta N_s=0\end{array} \right. .$$ 
Therefore, by independence and our representation
\begin{eqnarray}
&&E^{P^1}\left[\prod\limits _{u<s\le t}\left[1+\left(\frac{q_{Y_{s-}\rightarrow Y_{s}}\left(X_{s} \right)}{\overline q_{Y_{s-}\rightarrow Y_{s}}}-1\right)\Delta N_{s}\right]\Big|\mathcal F_u\right]\\\nonumber
&=&\!E^{P^1}\left[\prod\limits _{u<s\le t}\left[1+\left(\sum_j q_{Y_{s-}\rightarrow j}\left(X_{s} \right)-1\right)\Delta N_{s}\right]\Big|\mathcal F_u\right]=1
\end{eqnarray}
and $A^2 $ is a $P^1$-martingale.\\
(\emph{Use new m.p.}) (\ref{May30_22b}) is our intermediate m.p. 
We know that
\begin{eqnarray}\label{May20_22f}
g\left(Y_{t}\right)-g\left(Y_{0}\right)&=&\int _{0}^{t}\left[g\left(Y_{s}\right)-g\left(Y_{s-}\right)\right]dN_{s}\\\nonumber
&=&\int _{0}^{t}\left[g\left(\xi ^{Y_{s-}}_{N_{s-}}\right)-g\left(Y_{s-}\right)\right]dN_{s}\\\nonumber
&=&\int _{0}^{t}\left[g\left(\xi ^{Y_{s-}}_{N_{s-}}\right)-g\left(Y_{s-}\right)\right]\gamma_{Y_{s-}\rightarrow }\left(X_{s} \right)ds+\mathcal M_t
\end{eqnarray}
and $\mathcal M$ is a local $P^1$-martingale since it is the sum
of the following two terms, the first
\[
\int _{0}^{t}\left[\sum_{j}\overline q_{Y_{s-}\rightarrow j}\,g\left(j\right)-g\left(Y_{s-}\right)\right][dN_s-\gamma_{Y_{s-}\rightarrow }\left(X_{s} \right)ds]
\]
would clearly be a local martingale by (\ref{May30_22d}) and the second
\[
E^{P^1}\left[\int _{u}^{t}\left[g\left(\xi^{Y_{s-}}_{N_{s-}}\right)-\sum_{j}\overline q_{Y_{s-}\rightarrow j}g\left(j\right)\right][dN_s-\gamma_{Y_{s-}\rightarrow }\left(X_{s} \right)ds]\Big|\mathcal F_u\right]=0
\]
by the independence (from everything) and distribution of the $\{\xi_n^k\}$
under $P^1$.\\
Now, by (\ref{May30_22d},\ref{May20_22f})
\begin{eqnarray}\nonumber
&\!\!\!\!\!&\!\!\!\!\!\int _{0}^{t}\frac{q_{Y_{s-}\rightarrow Y_{s}}\left(X_{s} \right)}{\overline q_{Y_{s-}\rightarrow Y_{s}}}(g(Y_s)-g(Y_{s-}))d N_{s}-\!\int _{0}^{t}\frac{q_{Y_{s-}\rightarrow Y_{s}}\!\left(X_{s} \right)}{\overline q_{Y_{s-}\rightarrow Y_{s}}}\left[g\!\left(\xi ^{Y_{s-}}_{N_{s-}}\right)-g(Y_{s-})\right]\!\gamma_{Y_{s-}\rightarrow }\left(X_{s}\right)ds\\
&\!\!\!\!&\!\!\!\!=\widehat M^N_t\nonumber
\end{eqnarray}
is a $P^1$-local martingale.
Furthermore, if $\mathcal G_s=\sigma\{Y_v:v\le s\}\vee\mathcal F_u$ for $s>u$, then 
\begin{eqnarray}\nonumber
&\!\!\!\!&\!\!\!E\!\left[\int _{u}^{t}\frac{\ \gamma_{Y_{s-}\rightarrow \xi ^{Y_{s-}}_{N_{s-}}}(X_{s})}{\overline q_{Y_{s-}\rightarrow \xi ^{Y_{s-}}_{N_{s-}}}}\left[g\!\left(\xi ^{Y_{s-}}_{N_{s-}}\right)-g\!\left(Y_{s-}\right)\right]-\sum_j\left[g\left(j\right)-g\!\left(Y_{s-}\right)\right]\gamma_{Y_{s-}\rightarrow j}\left(X_{s} \right)ds\Big|\mathcal F_u\right]\\\nonumber
&\!\!\!\!=&\!\!\!\int _{u}^{t}\!E\!\left[E\!\left[\frac{\gamma_{Y_{s-}\rightarrow \xi ^{Y_{s-}}_{N_{s-}}}(X_{s})}{\overline q_{Y_{s-}\rightarrow \xi ^{Y_{s-}}_{N_{s-}}}}\left[g\!\left(\!\xi ^{Y_{s-}}_{N_{s-}}\right)-g\!\left(Y_{s-}\right)\right]-\!\sum_j\left[g(j)-g\!\left(Y_{s-}\right)\right]\gamma_{Y_{s-}\rightarrow j}\left(X_{s} \right)\!\Big|\mathcal G_{s-}\right]\!\Big|\mathcal F_u\right]\!ds\\
&\!\!\!\!=&\!0\nonumber
\end{eqnarray}
Consequently, 
\begin{equation}\label{Jun20_23f}
\widetilde M^N_t=\int _{0}^{t}\frac{q_{Y_{s-}\rightarrow Y_{s}}\left(X_{s} \right)}{\overline q_{Y_{s-}\rightarrow Y_{s}}}(g(Y_s)-g(Y_{s-}))d N_{s}-\!\int _{0}^{t}\!\sum_j\left[g(j)-g\!\left(Y_{s-}\right)\right]\gamma_{Y_{s-}\rightarrow j}\left(X_{s} \right)ds
\end{equation}
is a $P^1$-local martingale.\\
(\emph{Apply Girsanov-Meyer}) By (\ref{Jun20_23f}), Theorem III.20 of Protter and the fact that $A^2$ is pure jump, one has that
\begin{eqnarray}
\ \ \ \ \ \ \check M^N_t\label{May20_22e}
&\!\!\!=&\!\!\int _{0}^{t}\frac{q_{Y_{s-}\rightarrow Y_{s}}\left(X_{s} \right)}{\overline q_{Y_{s-}\rightarrow Y_{s}}}(g(Y_s)-g(Y_{s-})\!)d N_{s}-\!\int _{0}^{t}\!\sum_j\left[g(j)-g\!\left(Y_{s-}\right)\right]\gamma_{Y_{s-}\rightarrow j}\!\left(X_{s} \right)\!ds\!\\\nonumber
&\!\!\!-&\!\!\int_0^t\frac1{A^2_s}d[A^2,\int_0^{\cdot} \frac{q_{Y_{u-}\rightarrow Y_{u}}\left(X_{u} \right)}{\overline q_{Y_{u-}\rightarrow Y_{u}}}(g(Y_u)-g(Y_{u-}))d N_{u}]_s
\end{eqnarray}
is a local martingale under $P$.  
However, by (\ref{A2diffdef}) and (\ref{May20_22f})
\begin{eqnarray}\label{May20_22g}
&\!\!\!\!&\!\!\left[A^{2},\int \frac{q_{Y_{u-}\rightarrow Y_{u}}\left(X_{u} \right)}{\overline q_{Y_{u-}\rightarrow Y_{u}}}(g(Y_u)-g(Y_{u-}))d N_{u}\right]_{t}\\\nonumber
&\!\!\!\!=&\!\!\sum\limits _{0<s\le t}\left[A^2_{s-}\!\left(\frac{q_{Y_{s-}\rightarrow Y_{s}}\left(X_{s} \right)}{\overline q_{Y_{s-}\rightarrow Y_{s}}}-1\!\right)\frac{q_{Y_{s-}\rightarrow Y_{s}}\left(X_{s} \right)}{\overline q_{Y_{s-}\rightarrow Y_{s}}}(g(Y_s)-g(Y_{s-}))\right]\!\Delta N_{s}\ 
\end{eqnarray}
so by (\ref{May20_22b})
\begin{eqnarray}\label{May20_22g}
&\!\!\!\!&\!\!\int_0^t\frac1{A^2_s}d\left[A^{2},\int \frac{q_{Y_{u-}\rightarrow Y_{u}}\left(X_{u} \right)}{\overline q_{Y_{u-}\rightarrow Y_{u}}}(g(Y_u)-g(Y_{u-}))d N_{u}\right]_{s}\\\nonumber
&\!\!\!\!=&\!\!\int _{0}^{t}\frac{A^2_{s-}}{A^2_s}\left(\frac{q_{Y_{s-}\rightarrow Y_{s}}\left(X_{s} \right)}{\overline q_{Y_{s-}\rightarrow Y_{s}}}-1\!\right)\frac{q_{Y_{s-}\rightarrow Y_{s}}\left(X_{s} \right)}{\overline q_{Y_{s-}\rightarrow Y_{s}}}(g(Y_s)-g(Y_{s-}))d N_{s}\\
&\!\!\!\!=&\!\!\int _{0}^{t}\frac{\overline q_{Y_{s-}\rightarrow Y_{s}}}{q_{Y_{s-}\rightarrow Y_{s}}\left(X_{s} \right)}\left(\frac{q_{Y_{s-}\rightarrow Y_{s}}\left(X_{s} \right)}{\overline q_{Y_{s-}\rightarrow Y_{s}}}-1\!\right)\frac{q_{Y_{s-}\rightarrow Y_{s}}\left(X_{s} \right)}{\overline q_{Y_{s-}\rightarrow Y_{s}}}(g(Y_s)-g(Y_{s-}))d N_{s}.\nonumber
\end{eqnarray}
Substituting (\ref{May20_22g}) into (\ref{May20_22e}) and using bounded $g$, one has that
\begin{eqnarray}
\!\!\check M^N_t\label{newMcheck}
&\!\!=&\!\int _{0}^{t}(g(Y_s)-g(Y_{s-}))d N_{s}-\int _{0}^{t}\!\sum_j\left[g(j)-g\!\left(Y_{s-}\right)\right]\gamma_{Y_{s-}\rightarrow j}\left(X_{s} \right)ds\\ 
&\!\!=&\!g(Y_t)-g(Y_{0})-\int _{0}^{t}\!\sum_j\left[g(j)-g\!\left(Y_{s-}\right)\right]\gamma_{Y_{s-}\rightarrow j}\left(X_{s} \right)ds\nonumber
\end{eqnarray}
is a martingale under $P$.  \\
(\emph{Work on }$X$) A third application of the Girsanov-Meyer theorem with $[f(X),N]=0$ shows
\begin{equation}\label{Junef}
m^f_t=f\left(X_{t} \right)-\int _{0}^{t}Lf\left(X_{s}\right)ds
\end{equation}
is also a $P$ local martingale.  
One now obtains from boundedness, integration by parts, (\ref{newMcheck}), the fact $\left[m^{f},N\right]=0$ and (\ref{Junef}) that
\begin{eqnarray}
f\left(X_{t}\right)g\left(Y_{t}\right)-f\left(X_{0}\right)g\left(Y_{0}\right)
-\int _{0}^{t}g(Y_{s})Lf(X_{s})+f(X_{s})\mathbb L(X_s)g(Y_s)ds\ 
\end{eqnarray}
is a $P$ martingale for bounded $g$ and $f\in D(L)$.
\eproof

As mentioned previously, CTHMM is a special case of our CMOM model.
(C1-C3) are modified for the CTHMM model as follow:
\begin{itemize}
\item[(A1)] The observation state space $O$ is a finite or countable space.
\item[(A2)] $\sup\limits_{x\in E}\gamma(x)<\infty$.
\item[(A3)] $\gamma(x)>0$ for all $x\in E$.
\end{itemize}
Now, the following result is an immediate corollary of Theorem \ref{ImportanceTheorem}.
\begin{corollary}\label{CTHMMImportanceTheorem}
Suppose (A1, A2, A3) hold and $\left(X,Y\right)$ satisfies the (\ref{CTHMMStar}) m.p. starting from
some initial law $\mathcal L(X_0,Y_0)=\nu$ under $Q$.
Then, 
\begin{equation}\label{CTHMMAweight}
A_{t}=\exp\left(\int _{0}^{t}\overline\gamma-\gamma\left(X_{s} \right)ds\right)\prod\limits _{0<s\le t}\left[1+\left(\frac{\gamma(X_{s}) q_{Y_{s}}(X_{s})}{\overline\gamma\, \overline q_{Y_{s}}}-1\right)\Delta N_{s}\right]
\end{equation}
is a $\left\{\mathcal F^{Y}_{t}\right\}$-martingale under $Q$, where $N$ counts
the transitions of $Y$. 
Moreover, if we define $P$ via
$
\left.\frac{dP}{dQ}\right|_{\mathcal F_{t}}=A_{t},\ \ \forall t\ge 0,
$
then $(X,Y)$ satisfies the (\ref{CTHMMPlus}) m.p. starting from
$\mathcal L(X_0,Y_0)=\nu$ under $P$.  
\end{corollary}

\section{Simulation and Model Testing}

\subsection{Rejection Sampling}
Theorem \ref{ImportanceTheorem} was partially motivated by the desire 
to change the rates of Markov chains.  
We will not specify the exact m.p. domains for this application
for maximum generality.
We start with a (simple Markov chain say) m.p. under the reference 
probability measure $Q$ as the proposal process $Y$ satisfying:
\begin{equation}\label{SimpStar}
m^Q_t=g\left(Y_{t}\right)-\int _{0}^{t}\overline {\mathbb L}g\left(Y_{s}\right)ds
\end{equation}
is a martingale for all $g$.  
This can be expanded to the m.p.
\begin{equation}\label{Starmod}
M^Q_t=f\left({t} \right)g\left(Y_{t}\right)-\int _{0}^{t}g\left(Y_{s}\right)\frac{d}{ds}f\left({s} \right)ds-\int _{0}^{t}f\left(s \right)\overline {\mathbb L}g\left(Y_{s}\right)ds,
\end{equation}
is a martingale for all $f,g$ by integration by parts,  
which is just (\ref{Star}) with $X_s=s$.
(\ref{Starmod}) is then turned into target m.p. (\ref{Plus}) by Theorem
\ref{ImportanceTheorem} with $X_s=s$:
\begin{equation}\label{NPlus}
M^P_t=M^P_t(f,g)=f\left({t} \right)g\left(Y_{t}\right)-\!\int _{0}^{t}g\left(Y_{s}\right)\frac{d}{ds}f\left({s} \right)ds-\!\int _{0}^{t}f\left({s} \right)\mathbb L_{s} g\left(Y_{s}\right)ds\
\end{equation}
is a martingale for all $f,g$ under the new measure $P$.
The likelihood ratio martingale for this change is then given by (\ref{FirstA})
according to Theorem \ref{ImportanceTheorem}.
 
In summary, we
\begin{enumerate}
\item Simulate Markov chain $Y$ under the (simple) proposal process distribution.  We think of this as being done on the reference probability space
$\left(\Omega ,\mathcal F,Q\right)$.
\item Reweight the simulation by $A$ so that the combined effect is like it came from the target distribution with a different probability $P$.
\end{enumerate}
It is natural to wonder if there is some way to stay on the simulation space (with $Q$) and get rid of the weight.  
Notice sample $Y$ dependence of the likelihood ratio weight:
\begin{equation}\label{Ainspect}
A=A_{T}\left(Y\right)=\exp\left(\int _{0}^{T}\!\overline\gamma _{Y_{s}\rightarrow }-\gamma _{Y_{s}\rightarrow }(s)ds\!\right)\prod\limits _{0<s\le T}\!\left[1+\left(\!\frac{\gamma _{Y_{s-}\rightarrow Y_{s}}(s)}{\overline\gamma _{Y_{s-}\rightarrow Y_{s}}}-1\!\right)\Delta N_{s}\right]\!,\!\!\!
\end{equation}
where $N$ counts the transitions of $Y$, to convert proposal Markov chain simulations 
as target ones.
In fact, the weight tells us how good a proposal sample would be as a sample 
from the target m.p.
Now, von Neumann's Acceptance-Rejection algorithm:\vspace*{0.25cm}\\
(Step 1) Simulate $Y^Q$ with proposal m.p. and a $[0,C]$-Uniform $U$ independent of $Y^Q$.\vspace*{0.2cm}\\
(Step 2) If $U\le A_T\left(Y^{Q}\right)$, then accept by setting $Y=Y^Q$ and quitting the algorithm.  Otherwise, reject by returning to Step 1.
\vspace*{0.25cm}
\par\noindent
will allow us to create
a target sample without requiring the likelihood ratio weight.
We constrain $A$ to be bounded (for now) in order to show the algorithm works.
\vspace*{0.25cm}
\par\noindent
The following conditions will be imposed in our simulation result:
\begin{itemize}
\item[(C2')] $\sup\limits_{s\in [0,T],i\in O}\frac{\gamma_{i\rightarrow }(s)}{\overline\gamma_{i\rightarrow }}<\infty$.
\item[(C3')] There are no cemetery states, meaning $\gamma_{i\rightarrow }(s), \overline\gamma_{i\rightarrow }>0$ for all $i\in O$, $s\in [0,T]$.
\item[(Good $A$)] There is a $C>1$ such that $A=A_T\left(Y\right)\in \left[0,C\right]$ for all samples $Y$ and $E^{Q}\left[A\right]=1$.
\end{itemize}

\begin{proposition}
Suppose (C1, C2', C3', Good $A$) hold, $A$ is defined in (\ref{Ainspect}) and 
$\{(Y^Q_n,U_n)\}_{n=1}^\infty$ are the independent samples produced by Step 1
of the rejection algorithm. 
Then, the rejection algorithm output $Y$ has the target distribution on $Q$.
\end{proposition}
For clarity, the theorem says one can simulate on a computer with
rates $\{\overline \gamma_{i\rightarrow j}\}$ apply rejection resampling
and get a Markov chain with rates $\{\gamma_{i\rightarrow j}(s)\}$
at times $s$.
\proof
Let $f:\mathbb R^T\rightarrow \mathbb R$ be continuous and bounded. 
By the algorithm's acceptance criterion and independence one has (with $A^n=A_T(Y^Q_n)$) that
\begin{eqnarray}
E^Q[f(Y)]&=&\sum\limits _{n=1}^{\infty }E^{Q}\left[f\left(Y^{Q}_{n}\right)1_{U_{1}>A^{1}},...,1_{U_{n-1}>A^{n-1}},1_{U_{n}\le A^{n}}\right]\\\nonumber
&=&\sum\limits _{n=1}^{\infty }E^{Q}\left[f\left(Y^{Q}_{n}\right),1_{U_{n}\le A^{n}}\right]Q\left({U_{1}>A^{1}},...,{U_{n-1}>A^{n-1}}\right)\\\nonumber
&=&\sum\limits _{n=1}^{\infty }\frac{E^{Q}\left[f\left(Y^{Q}_{n}\right)A^{n} \right]}C Q\left({U_{1}}>A^{1}\right)\cdots Q\left({U_{n-1}}>A^{n-1}\right)\\\nonumber
&=&\sum\limits _{n=1}^{\infty }\frac{E^{Q}\left[f\left(Y^{Q}_{1}\right)A^{1} \right]}C Q\left({U_{1}}>A^{1}\right)^{n-1}\\\nonumber
&=&E^{Q}\left[f\left(Y^{Q}_{1}\right)A^{1} \right]\frac{1}{C Q\left({U_{1}}\le A^{1}\right)}\\\nonumber
&=&E^{P}\left[f\left(Y^{Q}_{1}\right) \right]\frac{1}{C Q\left({U_{1}}\le A^{1}\right)}.
\end{eqnarray}
Substituting $f=1$, one finds that $1=\frac{1}{CQ\left(U_{1}\le A^{1}\right)}$ so $Q\left(U_{1}\le A^{1}\right)=\frac{1}{C}$. $\square$

The acceptance rate decreases as $C$ increases.
Unfortunately, there is no absolute bound on the number of transitions $N$.
If $\left\{\gamma _{i\rightarrow j}\right\}$ are chosen so that $\frac{\gamma _{Y_{s-}\rightarrow Y_{s}}\left(s\right)}{\gamma _{Y_{s-}\rightarrow Y_{s}}}\le 1$ for all samples $Y$ and all $s$, then we can use the rejection method as above without stopping times since  
$$A=\exp\left(\int _{0}^{T}\overline\gamma _{Y_{s}\rightarrow }-\gamma _{Y_{s}\rightarrow }\left(s\right)ds\right)\prod\limits _{0<s\le T}\left[1+\left(\frac{\gamma _{Y_{s-}\rightarrow Y_{s}}\left(s\right)}{\overline\gamma _{Y_{s-}\rightarrow Y_{s}}}-1\right)\Delta N_{s}\right]$$
can be bounded.
(This can be either impossible or else force a very inefficient simulation.)
Otherwise, we use the algorithm multiple times with stopping times.

Let $W_n$ be the time of the $n^{th}$ jump of $N$, which is the $n^{th}$ transition of $Y$ and assume\\
(Bounded transitions) Suppose there is a $c>0$ such that $\left|\lambda _{i\rightarrow }\right|\le c$ and $\sup_{i,j}\left|\frac{\gamma _{i\rightarrow j}\left(s\right)}{\lambda _{i\rightarrow j}}\right|\le c$.

Now, the following lemma is trivial:
\begin{lemma}
Suppose (Bounded transitions) holds, $T>0$ and $n\in \mathbb{N}$.  Then, $A_{T\wedge W_{n}}$ satisfies (good $A$).
\end{lemma}
Indeed, it can be a good idea to simulate a fixed number of transitions at a time based upon the following formula.
$Y$ only changes when $N$ jumps and $N$ only jumps by $1$ and we can consider our explicit formula along the jump times $\left\{W_{n}\right\}$ of $N$.  In which case we get:
$$A_{W_{n}}=A_{W_{n-1}}\exp\left(\int ^{W_{n}}_{W_{n-1}}\left(\overline\gamma _{Y_{W_{n-1}}\rightarrow }-\gamma _{Y_{W_{n-1}}\rightarrow }\left(s\right)\right)ds\right)\left(\frac{\gamma _{Y_{W_{n-1}}\rightarrow Y_{W_{n}}}\left(W_{n}\right)}{\overline\gamma _{Y_{W_{n-1}}\rightarrow Y_{W_{n}}}}\right)$$
which shows how the explicit solution weight updates at jump times.  Here, we have that $W_{0}=0$ and $A_{0}=1$.
In this way, we can resample a few jump times at a time under
(Bounded transitions), provided we do not exceed some fixed time $T>0$.

\subsection{Monte Carlo Simulation}
Very often, we are not after just one sample but rather an ensemble of samples to form a distribution.  
We are really after some probabilities, expectations or conditional expectations, which we can approximate by independent proposal particles on the reference probability space $\left(\Omega ,\mathcal F,Q\right)$ weighted by likelihood martingale weights.  
In particular, suppose $\left\{Y^{m}\right\}_{m=1}^{M}$ are independent proposal particles on $\left(\Omega ,\mathcal F,Q\right)$ with rates $\left\{\overline\gamma _{i\rightarrow j}\right\}$ say and we weight the $m^{th}$ particle with likelihood weight
$$
A^{m}=A^{m}_{T}=\exp\left(\int _{0}^{T}\overline\gamma _{Y^{m}_{s}\rightarrow }-\gamma _{Y^{m}_{s}\rightarrow }\left(s\right)ds\right)\prod\limits _{0<s\le T}\left[1+\left(\frac{\gamma _{Y^{m}_{s-}\rightarrow Y^{m}_{s}}\left(s\right)}{\overline\gamma _{Y^{m}_{s-}\rightarrow Y^{m}_{s}}}-1\right)\Delta N^{m}_{s}\right],
$$
where $N^{m}$ counts $Y^m$'s jumps.  
Then, it follows by the strong law of large numbers that
$$\frac{1}{M}\sum\limits _{m=1}^{M}A^{m}_{T}f\left(Y^{m}_{s},s\in \left[0,T\right]\right)\rightarrow E^{Q}\left[A_{T}f\left(Y_{s},s\in \left[0,T\right]\right)\right]=E^{P}\left[f\left(Y_{s},s\in \left[0,T\right]\right)\right].$$
Hence, \emph{any} target distribution expectation can be estimated without (rejection) resampling.

\subsection{Particle Filtering and Model Selection}

Particle filters use Monte Carlo methods to approximate filters.
Branching particle filters are often among the best performers and have
the added advantages of readily providing unnormalized filters and Bayes
factor.
The idea of these particle filters is to start with the weighted
particle filter (\ref{WeightPart}) but then, as mentioned in the 
introduction, branch in an unbiased way so as to keep the weights
somewhat uniform while still keeping the number of particles 
relatively constant.
The following algorithm is an adaptation to our setting of the simple residual algorithm introduced and explained in 
\cite{Ko17a}.
The other algorithms given in \cite{Ko17a} could also be employed to improve
performance.

{\bf Algorithm Setting:}
Let $\{t_n\}_{n=1}^\infty$ be the random transition times of the observations and set $t_0=0$.  
Let $r\in \left(1,\infty \right)$ and $\left\{V_{n}^i\right\}_{n,i=1}^{\infty }$, $\left\{U_{n}^i\right\}_{n,i=1}^{\infty }$ be independent  $\left[-0.1,0.1\right]$-uniform, $[0,1]$-uniform random variables. 
$r$ is a resampling or branching parameter.  
$\left\{V_{n}^i\right\}$ are smoothing random variable's to ensure we can show convergence.  They can be tighter than $[-0.1,0.1]$. 

All particles evolve independently of each other between observations and interact weakly (i.e. through their empirical measure) at observation times.  

\vspace*{0.1cm}
\par\noindent
\textbf{Let}: $\left\{X_{0}^{i}\right\}_{i=1}^{N}$ be independent, $\mathcal L(X_{0}^{i})=\mu $; $N^{i}_{0}=A^{i}_{0}=1$ for $i=1,...,N$; $N_{n}=0$ for $n\in \mathbb{N}$.
\vspace*{0.2cm}
\par\noindent
\textbf{Repeat}: for $n=1,2,...$ do  
\vspace*{0.2cm}

  \textbf{Repeat}: for $i=1,...,N_{n-1}$ do
\begin{enumerate}
\item Evolve Particle Behavior Independently: Simulate each $X^i$ on $(t_{n-1},t_n]$ independently according to the signal's generator.  Call final point $\widetilde X^i_{t_{n}}$.

\item Get observation $Y_{t_{n}}$, which also gives the value of
$Y_s$ on $[t_n,t_{n+1})$.

\item Weight Particles: $\widehat A^{i}_{t_{n}}=A^{i}_{t_{n-1}}\exp\left(\int _{t_{n-1}}^{t_{n}}\overline\gamma _{Y_{s\rightarrow} }-\gamma_{Y_{s\rightarrow} }\left(X^i_{s} \right)ds\right)\frac{\gamma_{Y_{t_{n-1}}\rightarrow Y_{t_n}}\left(X^i_{t_n}\right)}{\overline\gamma_{Y_{t_{n-1}}\rightarrow Y_{t_n}}}$

\item Resample Decision: Let $\overline A_{t_{n}}=\frac{1}{N}\sum\limits _{i=1}^{N_{n}}\widehat A^{i}_{t_{n}}$.  If $\widehat A^{i}_{t_{n}}+V_{n}^i\notin \left(\frac{\overline A_{t_{n}}}r ,r \overline A_{t_{n}}\right)$ then the
Offspring Number, Weight are: $N^{i}_{n}=\left\lfloor \frac{\widehat A^{i}_{t_{n}}}{\overline A_{t_{n}}}\right\rfloor +1_{U^i_n\le \left\{\frac{\widehat A^{i}_{t_{n}}}{\overline A_{t_{n}}}\right\}} $, $\widetilde A^{i}_{t_{n}}=\overline A_{t_{n}}$; or else the
Offspring Number, Weight are: $N^{i}_{n}=1$, $\widetilde A^{i}_{t_{n}}=\widehat A^{i}_{t_{n}}$.

\item Resample: $A^{N_{n}+j}_{t_{n}}=\widetilde A^{i}_{t_{n}}$,  $X^{N_{n}+j}_{t_{n}}=\widetilde X^i_{t_{n}}$ for $j=1,...,N^{i}_{n}$

\item Add Offspring Number: $N_{n}=N_{n}+N_{n}^{i}$
\end{enumerate}

The branching particle filter approximations of $\sigma _{t_{n}}$ are then
$$S ^{N}_{t_n}\left(f\right)=\frac{1}{N}\sum\limits _{i=1}^{N_n}A^{i}_{t_n}f\left(X^{i}_{t_n}\right).$$
Notice compared to regular Monte Carlo we
only have one observable process $Y$ but many hidden processes
$\{X^i\}_{i=1}^N$, called particles.
Compared to the weighted particle filter, the particles are adjusted 
so they are more effective and a better approximation 
is achieved with the same number of particles.
Based upon \cite{Ko17b}, the following results are
expected.  
\begin{conjecture}\label{resampclt}
Under general regularity conditions, for any $n\in\mathbb N$,
the above Residual Branching particle filter satisfies:
\begin{description}
\item[Slln]
${ S}^{N}_{t_n}\Rightarrow \sigma _{t_n}$ (i.e.\ weak convergence) as
$N\rightarrow \infty$ \ a.s. [$Q^Y$];
\item[Mlln]
$\left|{ S}^{N}_{t_n}\left(f\right)-\sigma_{t_n}\left(f\right)\right|
\stackrel{N}{\ll}
N^{-\beta}$ \ a.s. [$Q^Y$] $\forall f\in \overline{C}(E)_+$, $0\le\beta<\frac12$.
\end{description}
\end{conjecture}
 
For related background in nonlinear filtering and sequential Monte Carlo, the 
reader is referred to the books \cite{Xiong08} and \cite{Chop20} as well as the vast literature, including
\cite{FKK72}, \cite{Zakai}, \cite{KuOc88}, \cite{KoLo08}, \cite{DMKoMi}, 
\cite{Chop04}, \cite{MaNe}, \cite{VLKuNe}, \cite{ElTo21} and their references and citations. 

\begin{remark}
Only one step changes in the case of CTHMM.
Step 3 becomes:\\
Weight Particles: $\widehat A^{i}_{t_{n}}=A^{i}_{t_{n-1}}\exp\left(\int _{t_{n-1}}^{t_{n}}\overline\gamma -\gamma\left(X^i_{s} \right)ds\right)
\frac{\gamma(X^i_{t_n}) q_{Y_{t_n}}(X^i_{t_n})}{\overline\gamma\, \overline q_{ Y_{t_n}}}$.\\
The rest of the algorithm is unchanged.
\end{remark}

\section{Filtering equations}

Our filtering equation approach uses the
{\em reference probability} $Q$ and related {\em unnormalized
filter} process $\sigma$. 
$P$ restricted to ${\mathcal
F}_t$ is absolutely continuous with respect to $Q$ restricted to
${\mathcal F}_t$ for $t\ge 0$ and the
observation process $Y$ is independent of the hidden state $X$ under $Q$. 
Let $E^Q[\cdot]$ denote expectation
with respect to $Q$, and consider
the additional regularity condition:
\begin{itemize}
\item[(U)] 
$X$ is a Markov chain with state space $E\subset \mathbb N$.
\end{itemize}
Some additional regularity is required to establish
uniqueness of (\ref{DMZsingf}).
We chose to restrict 
$X$ to be a Markov chain is 
immediately verifiable and already built into
the CTHMM model.

The observations $Y$ are much simpler under $Q$ so our strategy is to
work under the reference probability and first derive an equation for 
$\sigma_t$. 
Then,
apply It\^o's formula to obtain the equation for the desired
conditional distribution given by Bayes' formula
\[
\pi _t(f)=\sigma _t(f)/\sigma _t(1).
\]
Now, we can state our second main result, which is on the filtering equations.

\begin{theorem}\label{FiltEq}
Suppose (C1, C2, C3) hold, $\left(X,Y\right)$ satisfies the (\ref{Star}) martingale problem starting from
some initial law $\mathcal L(X_0,Y_0)=\nu$ under $Q$ and 
$$
\left.\frac{dP}{dQ}\right|_{\mathcal F_{t}}=A_{t},\ \ \forall t\ge 0,
$$
where $A$ is defined in (\ref{Aweight}). 
Then, $\sigma $, defined by $\sigma _{t}\left(B\right)=E^Q\left(A_t1_{X_{t}\in B}\big|\mathcal F^{Y}_{t}\right)$, solves (\ref{DMZWea}) and 
$\pi$, defined by $\pi _{t}\left(B\right)=P\left(X_{t}\in B\big|\mathcal F^{Y}_{t}\right)$, for Borel subsets $B$ of $E$, solves (\ref{FKK}).
Moreover, if (U) also holds, then $\sigma $ is the
unique strong $D_{{\mathcal M}_f(E)}[0,\infty)$-valued solution to 
(\ref{DMZsingf}).
\end{theorem}
\begin{remark}
We have the ideal situation of existence in the more general setting (\ref{DMZWea}) but uniqueness holding already in the narrow setting (\ref{DMZsingf}) for $\sigma$.
\end{remark}
\begin{remark}
In the case that $f$ only depends upon $X$, one has that
\begin{eqnarray}
\pi_t(f)&\!=&\!\pi_0 (f)
+\int_0^t\pi_s (Lf) ds 
- \int_0^t\pi_s(f \gamma_{Y_{s}\rightarrow })
-\pi_s (f)\pi_s (\gamma_{Y_{s}\rightarrow })ds\\\nonumber
&\!+ &\!\int_0^t 
\frac {\pi_{s-} \left(
\gamma_{Y_{s-}\rightarrow Y_s}f\right)
-\pi_{s-} \left(\gamma_{Y_{s-}\rightarrow Y_s} \right)\pi_{s-}(f)}
{\pi_{s-}\left(\gamma_{Y_{s-}\rightarrow Y_s}\right)}
dN_s
\end{eqnarray}
for all $f\in\overline C(E)$.
\end{remark}

{\bf Proof $\sigma$ satisfies (\ref{DMZWea}).}
One notes by $Q$-independence, our representation $Y_s=\xi_{N_{s-}}^{Y_{s-}}$ and integration by parts that
\begin{eqnarray}\label{May10_22a}
&&f(X_t,Y_t)-f(X_0,Y_0)\\\nonumber
&\!\!=&\!\!\int_0^t\! Lf(X_s,Y_s)ds+m_t(f)+
\!\int_0^t [f(X_{s-},\xi_{N_{s-}}^{Y_{s-}})-f(X_{s-},Y_{s-})]dN_s
\end{eqnarray}
so by (\ref{May10_22a}) and (\ref{May10_22b}) (with $\alpha=A$) as well as independence
\begin{equation}\label{BieqnQ}
\lbrack f (X,Y) ,A\rbrack_t=\int_{0}^t A_{s-}[f(X_{s-},\xi_{N_{s-}}^{Y_{s-}})-f(X_{s-},Y_{s-})]
\left[\frac{\gamma_{Y_{s-}\rightarrow Y_s}\left(X_{s} \right)}{\overline \gamma_{Y_{s-}\rightarrow Y_s }}-1\right]
dN_s.\ \ 
\end{equation}
($L$ only operates on the first variable of $f$ in (\ref{May10_22a}) and below.)
Utilizing integration by parts, (\ref{May10_22a}), (\ref{May10_22b}) and
(\ref{BieqnQ}), one finds that 
\begin{eqnarray}\label{IntbyParts}
&\!&\!f(X_t,Y_t)A_t-f(X_0,Y_0)\\\nonumber
&\!=&\!\int_0^tA_{s-}
df(X_s,Y_s)+\int_0^tf(X_{s-},Y_{s-})dA_s+[f(X,Y),A]_t\\\nonumber
&\!\!\!=&\!\!\int_0^tA_{s} Lf(X_s,Y_s)ds+\int_0^tA_{s-} dm_s(f)
\\\nonumber
&\!\!\!+&\!\!\int_0^tA_{s-}f(X_{s-},Y_{s-})\left(\frac {\gamma_{Y_{s-}\rightarrow Y_s}\left(X_{s} \right)}{\overline\gamma_{Y_{s-}\rightarrow Y_s}}-1\right)dN_s\\\nonumber
&\!\!\!+&\!\!\int_0^tA_{s}f(X_{s},Y_{s})(\overline\gamma_{Y_{s}\rightarrow }-{\gamma_{Y_{s}\rightarrow }\left(X_{s} \right)})
ds\\\nonumber
&\!\!\!+&\!\!\int_{0}^t A_{s-}[f(X_{s-},Y_{s})-f(X_{s-},Y_{s-})]
\frac{\gamma_{Y_{s-}\rightarrow Y_s}\left(X_{s} \right)}{\overline \gamma_{Y_{s-}\rightarrow Y_s }}
dN_s\\\nonumber
&\!\!\!=&\!\!\int_0^tA_{s} Lf(X_s,Y_s)ds+\int_0^tA_{s-} dm_s(f)
+\int_0^tA_{s}f(X_{s},Y_{s})(\overline\gamma_{Y_{s}\rightarrow }-{\gamma_{Y_{s}\rightarrow }\left(X_{s} \right)})
ds\\
\nonumber
&\!\!\!+&\!\!\int_{0}^t A_{s-}\left[f(X_{s-},Y_{s})\frac{\gamma_{Y_{s-}\rightarrow Y_s}\left(X_{s} \right)}{\overline \gamma_{Y_{s-}\rightarrow Y_s }}-f(X_{s-},Y_{s-})\right]dN_s.
\end{eqnarray}
Next, we show $ E^Q[\int_0^tA_{s-} dm_s(f)|{\mathcal F}^Y_t]=0$.
For each $n\in\mathbb N$, let $t^n_0=0$ and $\left\{t^n_{i}\right\} _{i=1}^{\infty }$ be a 
refining partition of stopping times that include the transition
times of $Y$ such that 
$$
A^n_{s-}\doteq 1_{\{0\}}(s)+\sum_{i=0}^n A_{t_i^n}1_{(t_i^n,t_{i+1}^n]}(s)
$$
satisfies
$$
\sup_{0\le s\le t}|A_{s-}-A^n_{s-}|\rightarrow 0\ \ \mbox{in probability for any }t>0.
$$
Then, $\int_0^tA^n_{s-} dm_s(f)\rightarrow\int_0^tA_{s-} dm_s(f)$ in probability
and 
\begin{equation}\label{PriorOne}
E^Q|\int_0^tA^n_{s-} dm_s(f)-\int_0^tA_{s-} dm_s(f)|\rightarrow0
\end{equation}
by the boundedness of $m(f)$ and Condition (C1).
Moreover, it follows by the tower property, independence
and Doob's Optional Stopping that
\begin{eqnarray}\label{DoobJun19}
\ \ \ \ E^Q\!\left[\int_0^tA^n_{s-} dm_s(f)\big|{\mathcal F}^Y_t\right]
&\!=&\!
\sum_i E^Q\!\left[E^Q\!\left[A_{t_i^n}(m_{t_{i+1}^n}(f)-m_{t_{i}^n}(f))\big|{\mathcal F}^X_{t_{i}^n}\vee{\mathcal F}^Y_t\right]\!\big|{\mathcal F}^Y_t\right]\!\!\!\\
\nonumber
&\!=&\!
\sum_i E^Q\!\left[A_{t_i^n}E^Q\!\left[(m_{t_{i+1}^n}(f)-m_{t_{i}^n}(f))\big|{\mathcal F}^X_{t_{i}^n}\right]\!\big|{\mathcal F}^Y_t\right]\\
\nonumber
&\!=&0\ \ a.s.
\end{eqnarray}
Thus, it follows by (\ref{PriorOne}), (\ref{DoobJun19}) and Jensen's inequality that
\begin{eqnarray}\label{CombJun19}
\ \ \ \ E^Q\bigg|E^Q\!\left[\int_0^tA_{s-} dm_s(f)\big|{\mathcal F}^Y_t\!\right]\!\bigg|
&\!\!=&\!\!\!\lim_{n\rightarrow\infty}\!E^Q\!\bigg|E^Q\!\left[\int_0^t\!A_{s-} dm_s(f)-\int_0^t\!A^n_{s-} dm_s(f)\big|{\mathcal F}^Y_t\!\right]\!\bigg|\!\!\!\!\\
\nonumber
&\!\!\le&\!\!\lim_{n\rightarrow\infty}E^Q\left[\left|\int_0^tA_{s-} dm_s(f)-\int_0^tA^n_{s-} dm_s(f)\right|\right] \\
\nonumber
&\!\!=&\!0.\ \
\end{eqnarray}
Letting $E^*$ denote $Q$-expectation
with respect to $X$ only, and setting
\[
\sigma_t(f)\equiv E^Q[f(X_{t},Y_t)A_t|{\mathcal F}^Y_t]=E^{*}[f(X_{t},Y_t)L_t],
\]
we find by (\ref{IntbyParts}) that the Zakai-type equation 
for $\sigma_t (f)$ becomes
\begin{eqnarray*}
\!\sigma_t (f(\cdot,Y_t))&\!=&\!\sigma_0 (f(\cdot,Y_0))+\int_0^t\sigma_s (Lf(\cdot,Y_s)) ds+\int_0^t\sigma_s (f(\cdot,Y_{s})(\overline\gamma_{Y_{s}\rightarrow }-{\gamma_{Y_{s}\rightarrow }\left(\cdot \right)})
ds \\
&\!&\!+\int_{0}^t 
\sigma_{s-}\left(\left[f(\cdot ,Y_s)\frac{\gamma_{Y_{s-}\rightarrow Y_s}\left(\cdot  \right)}{\overline \gamma_{Y_{s-}\rightarrow Y_s}}-f(\cdot ,Y_{s-})\right]
\right)dN_s.
\end{eqnarray*} 

{\bf Proof of (\ref{DMZsingf}) uniqueness}: 
For ease of notation, we take $E=\mathbb N$ and $t_{0}=0$.
Let $\left\{t_{l}\right\} _{l=1}^{\infty }$ be the random transition times 
(in order) for $Y$ and
$Lf(i)=\sum\limits_{j\ne i}\lambda_{i\rightarrow j}[f(j)-f(i)]$ be
$X$'s generator.
Then, the adjoint operator $L^* p(j)$ satisfies
$$
L^*p\left(j\right)=\sum\limits _{i}\left[\lambda _{i\rightarrow j}p\left(i\right)-\lambda _{j\rightarrow i}\ p\left(j\right)\right]=[L^*]p (j),
\quad
\left[L^*\right]=\left[\begin{array}{cccc}-\lambda _{1\rightarrow }&\lambda _{2\rightarrow 1}&\lambda _{3\rightarrow 1}&\cdots\\
\lambda _{1\rightarrow 2}&-\lambda _{2\rightarrow }&\lambda _{3\rightarrow 2}&\cdots\\
\lambda _{1\rightarrow 3}&\lambda _{2\rightarrow 3}&-\lambda _{3\rightarrow }&\cdots\\
\vdots&\vdots&\vdots&\ddots
\end{array}\right].
$$
Now, if we let $\sigma^i_t=\sigma_t(\delta_i)$ for $i=1,2,...$, then we will 
discover that (\ref{DMZsingf}) gives us the closed system of linear differential
equations
parameterized by the observations: 
\begin{eqnarray}
\!\!d\!\left[\begin{array}{c}\sigma^1_t\\\sigma^2_t\\\vdots\end{array}\!\right] &\!=&\label{DMZdirect}\!
\left[\begin{array}{ccc}\overline\gamma_{Y_{t}\rightarrow }-{\gamma_{Y_{t}\rightarrow }(1)}-\lambda _{1\rightarrow }&\lambda _{2\rightarrow 1}&\cdots\\
\lambda _{1\rightarrow 2}&\overline\gamma_{Y_{t}\rightarrow }-{\gamma_{Y_{t}\rightarrow }(2 )}-\lambda _{2\rightarrow }&\cdots\\
\vdots&\vdots&\ddots\end{array}\right]\!\!\left[\begin{array}{c}\sigma^1_t\\\sigma^2_t\\\vdots\end{array}\right]\! dt\ \\
&+&\nonumber
\left[\begin{array}{ccc}\frac{{\gamma_{Y_{t-}\rightarrow Y_t}\left(1 \right)}}{\overline\gamma_{Y_{t-}\rightarrow  Y_t}}-1&0&\cdots\\
0&\frac{{\gamma_{Y_{t-}\rightarrow Y_t}\left(2 \right)}}{\overline\gamma_{Y_{t-}\rightarrow  Y_t}}-1&\cdots\\
\vdots&\vdots&\ddots\end{array}\right]\!\left[\begin{array}{c}\sigma^1_{t-}\\\sigma^2_{t-}\\\vdots\end{array}\right]\!dN_t.
\end{eqnarray}
Using the \emph{mild solution} and
\emph{Trotter product}, we can write the solution explicitly between observation times.
Let $P_t(i\rightarrow j),\ i,j\in \{1,2,...,m\}$ be the transition function for the hidden Markov chain $X$.
Then, the unnormalized filter is
$\sigma_t(\cdot)=\sum\limits_{i\in E}\sigma^i_t\delta_i(\cdot)$, where \begin{eqnarray}
\left[\begin{array}{c}\sigma^1_t\\\sigma^2_t\\\vdots\\\end{array}\right] &\!=&\!\left[T^n_{t-t_{n-1}}\right]
\label{Step1}
\left[\begin{array}{c}\sigma^1_{t_{n-1}}\\ \sigma^2_{t_{n-1}}\\
\vdots\end{array}\!\right]
\end{eqnarray}
for all $t\in [t_{n-1},t_n)$ with $T^n_t=\lim\limits_{N\rightarrow\infty}\left[S^n_{\frac{t-t_{n-1}}N}\right]^N$
and
Trotter product factor
\begin{eqnarray}
\ \ \ \ \ S^n_{t}\!&=\!&\!\!\label{Trotterinf}
\left[\!\begin{array}{ccc}
P_{t}(1\rightarrow 1)&P_{t}(2\rightarrow 1)&\!\cdots\\
P_{t}(1\rightarrow 2)&P_{t}(2\rightarrow 2)&\!\cdots\\
\vdots&\vdots&\!\ddots\end{array}\!\right]\!\!\!
\left[\!\begin{array}{ccc}e^{t(\overline\gamma_{Y_{t_{n-1}}\rightarrow }\!-{\gamma_{Y_{t_{n-1}}\rightarrow }(1 )})}\!&0&\cdots\\0&\!e^{t(\overline\gamma_{Y_{t_{n-1}}\rightarrow }\!-{\gamma_{Y_{t_{n-1}}\rightarrow }(2 )})}\!&\cdots\\\vdots&\vdots&\ddots\end{array}\!\right]\!.\!\!\!\!
\end{eqnarray}
Each $T^n$ behaves as a semi-group on
$[0,t_n-t_{n-1})$. 
Then,
\begin{equation}
\!\!\!\left[\begin{array}{l}\sigma^1_{t_n}\\\sigma^2_{t_n}\\\vdots\end{array}\!\right]
\label{Step2inf}
=\left[\begin{array}{c}\frac{\gamma_{Y_{t_{n-1}}\rightarrow Y_{t_{n}}}(1 )}{\overline\gamma_{Y_{t_{n-1}}\rightarrow Y_{t_{n}}}}\sigma^1_{t_{n}-}\\\frac{\gamma_{Y_{t_{n-1}}\rightarrow Y_{t_{n}}}(2 )}{\overline\gamma_{Y_{t_{n-1}}\rightarrow Y_{t_{n}}}} \sigma^2_{t_{n}-}\\\vdots\end{array}\!\right],
\end{equation}
and the equations start at $\left[\begin{array}{c}\sigma^1_{t_0}\\\sigma^2_{t_0}\\\vdots\end{array}\!\right]=\left[\begin{array}{c}P(X_0=1)\\P(X_0=2)\\\vdots\end{array}\!\right]$.

Now, suppose $\sigma,\widehat \sigma $ are $D_{\mathcal{M}_{f}(E)}[0,\infty )$-valued solutions to (\ref{DMZsingf}) and $\sigma =\widehat \sigma$ 
on $\left[0,t_{n-1}\right] $. 
Then, 
\[
\left[\begin{array}{c}\sigma^1_{t}\\\sigma^2_{t}\\\vdots\end{array}\right]=T^n_{t-t_{n-1}}\left[\begin{array}{c}\sigma^1_{t_{n-1}}\\\sigma^2_{t_{n-1}}\\\vdots\end{array}\right]=T^n_{t-t_{n-1}}\left[\begin{array}{c}\widehat\sigma^1_{t_{n-1}}\\\widehat\sigma^2_{t_{n-1}}\\\vdots\end{array}\right]=\left[\begin{array}{c}\widehat\sigma^1_{t}\\\widehat\sigma^2_{t}\\\vdots\end{array}\right]
\]
for $t\in [t _{n-1},t _{n})$ so $\sigma_t(\cdot)=\sum\limits_{i\in E}\sigma^i_t\delta_i(\cdot)=\sum\limits_{i\in E}\widehat\sigma^i_t\delta_i(\cdot)=
\widehat\sigma_t(\cdot)$ and uniqueness holds on $[0,t _n)$. 
Finally, (\ref{Step2inf}) yields
\begin{equation}
\!\!\!\left[\begin{array}{c}\sigma^1_{t_n}\\\sigma^2_{t_n}\\\vdots\end{array}\!\right]
\label{Step2infa}
=\left[\begin{array}{c}\frac{\gamma_{Y_{t_{n-1}}\rightarrow Y_{t_{n}}}(1 )}{\overline\gamma_{Y_{t_{n-1}}\rightarrow Y_{t_{n}}}}\sigma^1_{t_{n}-}\\\frac{\gamma_{Y_{t_{n-1}}\rightarrow Y_{t_{n}}}(2 )}{\overline\gamma_{Y_{t_{n-1}}\rightarrow Y_{t_{n}}}} \sigma^2_{t_{n}-}\\\vdots\end{array}\!\right]
=\left[\begin{array}{c}\frac{\gamma_{Y_{t_{n-1}}\rightarrow Y_{t_{n}}}(1 )}{\overline\gamma_{Y_{t_{n-1}}\rightarrow Y_{t_{n}}}}
\widehat\sigma^1_{t_{n}-}\\\frac{\gamma_{Y_{t_{n-1}}\rightarrow Y_{t_{n}}}(2 )}{\overline\gamma_{Y_{t_{n-1}}\rightarrow Y_{t_{n}}}} 
\widehat\sigma^2_{t_{n}-}\\\vdots\end{array}\!\right]=\left[\begin{array}{c}
\widehat\sigma^1_{t_n}\\
\widehat\sigma^2_{t_n}\\\vdots\end{array}\!\right],
\end{equation}
so (\ref{DMZsingf}) (strong) uniqueness holds on $[0,t _{n}]$
and on $[0,\infty)$ by induction.

{\bf Proof of (\ref{FKK})}: 
It follows by (\ref{DMZWea}) that
\begin{eqnarray}\label{May9a}
\ \sigma_s(f(\cdot,Y_s))&\!=&\!\sigma_{s-}(f(\cdot,Y_{s-}))+
\sigma_{s-}\!\left(\left[f(\cdot ,Y_s)\frac{\gamma_{Y_{s-}\rightarrow Y_s}\left(\cdot  \right)}{\overline \gamma_{Y_{s-}\rightarrow Y_s}}-f(\cdot ,Y_{s-})\right]
\!\right)\!\Delta N_s
\end{eqnarray}
so
\begin{eqnarray}\label{May9b}
\frac{\sigma_{s}(1)}{\sigma_{s-}(1)}=1+
\pi_{s-}\left(\frac{\gamma_{Y_{s-}\rightarrow Y_s}\left(\cdot  \right)}{\overline \gamma_{Y_{s-}\rightarrow Y_s}}
-1\right) 
\Delta N_s\ .
\end{eqnarray}
Next, recalling $\pi _t(f)=
\frac {\sigma _t(f)}{\sigma _t(1)}$ and using (\ref{May9a}) twice then (\ref{May9b}),
one has that
\begin{eqnarray}
&&\displaystyle \pi_s (f(\cdot,Y_s))-\pi_{s-}(f(\cdot,Y_{s-}))\label{pistarstar}\\\nonumber
&=& \frac{\sigma_s(f(\cdot,Y_s))-\sigma_{s-}(f(\cdot,Y_{s-}))-(\sigma_s(1)-\sigma_{s-}(1))\pi_{s-}(f(\cdot,Y_{s-}))}{\sigma_s(1)}\\\nonumber
&=&
\frac {\sigma_{s-} \left(
\gamma_{Y_{s-}\rightarrow Y_s}\left(\cdot  \right)f(\cdot,Y_s)\right)
-\sigma_{s-} \left(\gamma_{Y_{s-}\rightarrow Y_s}\left(\cdot  \right) \right)\pi_{s-}(f(\cdot,Y_{s-}))}
{\sigma_{s}(1)\overline \gamma_{Y_{s-}\rightarrow Y_s} }\Delta N_s\\\nonumber
&=& 
\frac {\pi_{s-} \left(
\gamma_{Y_{s-}\rightarrow Y_s}\left(\cdot  \right)f(\cdot,Y_s)\right)
-\pi_{s-} \left(\gamma_{Y_{s-}\rightarrow Y_s}\left(\cdot  \right) \right)\pi_{s-}(f(\cdot,Y_{s-}))}
{\pi_{s-}\left(\frac{\gamma_{Y_{s-}\rightarrow Y_s}\left(\cdot  \right)}{\overline \gamma_{Y_{s-}\rightarrow Y_s}}
\right)\overline \gamma_{Y_{s-}\rightarrow Y_s} }\Delta N_s\\\nonumber
&=& \frac {\pi_{s-} \left(
\gamma_{Y_{s-}\rightarrow Y_s}\left(\cdot  \right)f(\cdot,Y_s)\right)
-\pi_{s-} \left(\gamma_{Y_{s-}\rightarrow Y_s}\left(\cdot  \right) \right)\pi_{s-}(f(\cdot,Y_{s-}))}
{\pi_{s-}\left(\gamma_{Y_{s-}\rightarrow Y_s}\left(\cdot  \right)\right)}
\Delta N_s.
\end{eqnarray}
Ito's formula on $\pi_t(f(\cdot,Y_t))=\frac{\sigma_t(f(\cdot,Y_t))}{\sigma_t(1)}$ gives
\begin{equation}\label{ItoSub}
d\pi_t(f(\cdot,Y_t))=\frac{d\sigma_t^c(f(\cdot,Y_t))}{\sigma_t(1)}-\pi_t(f(\cdot,Y_t))\frac{d\sigma_t^c(1)}{\sigma_t(1)}
+\pi_t(f(\cdot,Y_t))-\pi_{t-}(f(\cdot,Y_t)),
\end{equation}
where $c$ indicates the continuous part. 
Hence, using (\ref{DMZWea}), (\ref{ItoSub}) and (\ref{pistarstar}), 
one has
\begin{eqnarray*}
\pi_t(f(\cdot,Y_t))&\!=&\!\pi_0 (f(\cdot,Y_0))
+\int_0^t\pi_s (Lf(\cdot,Y_{s})) ds \\
&\!- &\!\int_0^t\pi_s(f(\cdot,Y_{s}) (\gamma_{Y_{s}\rightarrow }\left(\cdot  \right)
-\overline\gamma_{Y_{s}\rightarrow }))
-\pi_s (f(\cdot,Y_{s}))\pi_s (\gamma_{Y_{s}\rightarrow }\left(\cdot  \right)
-\overline\gamma_{Y_{s}\rightarrow })ds\\
&\!+ &\!\int_0^t 
\frac {\pi_{s-} \left(
\gamma_{Y_{s-}\rightarrow Y_s}\left(\cdot  \right)f(\cdot,Y_s)\right)
-\pi_{s-} \left(\gamma_{Y_{s-}\rightarrow Y_s}\left(\cdot  \right) \right)\pi_{s-}(f(\cdot,Y_{s-}))}
{\pi_{s-}\left(\gamma_{Y_{s-}\rightarrow Y_s}\left(\cdot  \right)\right)}
dN_s
\ . \hfill \qquad\ \  \Box
\end{eqnarray*}
For completeness, we state the corresponding filtering equations
result for the popular CTHMM, which is an immediate corollary.
However, to do this we have to first give the CTHMM versions of
(\ref{FKK}) and (\ref{DMZWea}), which are:
\begin{eqnarray}\label{CTHMMDMZWea}
\ \ \ \sigma_t (f(\cdot,Y_t))
&=&\sigma_0 (f(\cdot,Y_0))+\int_0^t\!\sigma_s (Lf(\cdot,Y_s)) ds+\int_0^t\!\sigma_s (f(\cdot,Y_{s})(\overline\gamma-{\gamma\left(\cdot \right)})
ds\\\nonumber
&&+\int_{0}^t 
\sigma_{s-}\left(\left[f(\cdot ,Y_s)\frac{\gamma(\cdot  ) q_{Y_s}(\cdot  )}{\overline \gamma \,\overline q_{ Y_s}}-f(\cdot ,Y_{s-})\right]
\right)dN_s
\end{eqnarray} 
for all $f\in\hat D_L$ and
\begin{eqnarray}\label{CTHMMFKK}
\ \ \ \ \ \pi_t(f(\cdot,Y_t))&\!=&\pi_0 (f(\cdot,Y_0))
+\int_0^t\pi_s (Lf(\cdot,Y_{s})) ds \\\nonumber
&\!- &\!\int_0^t\!\pi_s(f(\cdot,Y_{s}) \gamma\left(\cdot  \right)
)
-\pi_s (f(\cdot,Y_{s}))\pi_s (\gamma\left(\cdot  \right)
)ds\\\nonumber
&\!+ &\!\!\int_0^t 
\frac {\pi_{s-} \left(
\gamma\left(\cdot  \right) q_{ Y_s}\left(\cdot  \right)f(\cdot,Y_s)\right)
-\pi_{s-} \!\left(\gamma\left(\cdot  \right) q_{ Y_s}\left(\cdot  \right) \right)\pi_{s-}(f(\cdot,Y_{s-}))}
{\pi_{s-}\left(\gamma\left(\cdot  \right) q_{ Y_s}\left(\cdot  \right)\right)}
dN_s,
\end{eqnarray}
for all $f\in\hat D_L$ subject to $\pi_0=\mathcal L(X_0)$.
Now, the CTHMM corollary is:
\begin{corollary}\label{CTHMMFiltEq}
Suppose (A1, A2, A3) hold, $\left(X,Y\right)$ satisfies the (\ref{CTHMMStar}) martingale problem starting from
some initial law $\mathcal L(X_0,Y_0)=\nu$ under $Q$ and 
$
\left.\frac{dP}{dQ}\right|_{\mathcal F_{t}}=A_{t},\ \ \forall t\ge 0
$,
where $A$ is defined in (\ref{CTHMMAweight}). 
Then, $\sigma $ solves (\ref{CTHMMDMZWea}) and $\pi$ solves (\ref{CTHMMFKK}).
Moreover, if (U) also holds, then $\sigma $ is the
unique strong $D_{{\mathcal M}_f(E)}[0,\infty)$-valued solution to 
(\ref{CTHMMDMZ}).
\end{corollary}

\subsection{Direct Solution}

An effective computer workable solution to many real filtering problems can 
constructed from the DMZ equation based on uniqueness technique above.
Consider the case where $X$ is (or has been approximate by) a Markov chain on a 
finite space $E=\{1,2,...,m\}$ with generator
$Lf\left(i\right)=\sum\limits _{j\ne i}\lambda _{i\rightarrow j}
\left[f\left(j\right)-f\left(i\right)\right]$ for $i\in\{1,2,...,m\}$.
Now, if $\sigma^i_t=\sigma_t(\delta_i)$ for $i=1,2,...,m$, then (\ref{DMZsingf}) gives us the system of equations: 
\begin{eqnarray}
\!\!d\!\left[\!\begin{array}{l}\sigma^1_t\\\sigma^2_t\\\vdots\\\sigma^m_t\end{array}\!\!\right] &\!\!\!=&\!\!\!
\left[\!\begin{array}{cccc}\overline\gamma_{Y_{t}\rightarrow }\!\!-{\gamma_{Y_{t}\rightarrow }(1)}-\lambda _{1\rightarrow }&\lambda _{2\rightarrow 1}&\!\cdots&\lambda _{m\rightarrow 1}\\
\lambda _{1\rightarrow 2}&\overline\gamma_{Y_{t}\rightarrow }\!-{\gamma_{Y_{t}\rightarrow }(2 )}-\lambda _{2\rightarrow }\!&\!\cdots&\lambda _{m\rightarrow 2}\\
\vdots&\vdots&\!\ddots&\vdots\\
\!\lambda _{1\rightarrow m}&\lambda _{2\rightarrow m}&\!\cdots&\overline\gamma_{Y_{t}\rightarrow }\!-{\gamma_{Y_{t}\rightarrow }(m )}-\lambda _{m\rightarrow }\!\end{array}\right]\!\!\!\left[\!\begin{array}{l}\sigma^1_t\\\sigma^2_t\\\vdots\\\sigma^m_t\end{array}\!\!\right]\! dt\ 
\nonumber\\\nonumber
&\!\!\!+&\!\!\!\label{DMZdirect}
\left[\begin{array}{cccc}\frac{{\gamma_{Y_{t-}\rightarrow Y_t}\left(1 \right)}}{\overline\gamma_{Y_{t-}\rightarrow  Y_t}}-1&0&\cdots&0\\
0&\frac{{\gamma_{Y_{t-}\rightarrow Y_t}\left(2 \right)}}{\overline\gamma_{Y_{t-}\rightarrow  Y_t}}-1&\cdots&0\\
\vdots&\vdots&\ddots&\vdots\\
0&0&&\frac{{\gamma_{Y_{t-}\rightarrow Y_t}\left(m \right)}}{\overline\gamma_{Y_{t-}\rightarrow  Y_t}}-1\end{array}\right]\!\left[\begin{array}{l}\sigma^1_{t-}\\\sigma^2_{t-}\\\vdots\\\sigma^m_{t-}\end{array}\!\right]\!dN_t.
\end{eqnarray}
Let $\{t_n\}_{n=1}^\infty$ be the random transition times of the observations, $t_0=0$,
$P_t(i\rightarrow j),\ i,j\in \{1,2,...,m\}$ be the transition function for the hidden-state Markov chain $X$
and
\begin{eqnarray}
S^n_{t}&=\!&\label{Trotter}
\left[\!\begin{array}{cccc}
P_{t}(1\rightarrow 1)&P_{t}(2\rightarrow 1)&\cdots&P_{t}(m\rightarrow 1)\\
P_{t}(1\rightarrow 2)&P_{t}(2\rightarrow 2)&\cdots&P_{t}(m\rightarrow 2)\\
\vdots&\vdots&\ddots&\vdots\\
P_{t}(1\rightarrow m)&P_{t}(2\rightarrow m)&\cdots&P_{t}(m\rightarrow m)\end{array}\right]\!\!
\\\nonumber
&*&\!\left[\begin{array}{cccc}e^{t(\overline\gamma_{Y_{t_{n-1}}\rightarrow }-{\gamma_{Y_{t_{n-1}}\rightarrow }(1 )})}\!&0&\cdots&0\\0&e^{t(\overline\gamma_{Y_{t_{n-1}}\rightarrow }-{\gamma_{Y_{t_{n-1}}\rightarrow }(2 )})}\!\!&\cdots&0\\\vdots&\vdots&\ddots&\vdots\\0&0&\cdots&\!e^{t(\overline\gamma_{Y_{t_{n-1}}\rightarrow }\!-{\gamma_{Y_{t_{n-1}}\rightarrow }(m )})}\end{array}\!\right]
\end{eqnarray}
for some large $N$ (Trotter product semi-group approximation).
Then, the solution is found recursively by i) 
weighted evolution
\begin{eqnarray}
\left[\begin{array}{c}\sigma^1_t\\\sigma^2_t\\\vdots\\\sigma^m_t\end{array}\!\right] &=&\left[S^n_{\frac{t-t_{n-1}}N}\right]^N
\label{Step1}
\left[\begin{array}{c}\sigma^1_{t_{n-1}}\\ \sigma^2_{t_{n-1}}\\
\vdots\\
\sigma^m_{t_{n-1}}\end{array}\right]
\end{eqnarray}
for all $t\in [t_{n-1},t_n)$ and ii) observation outcome update
\begin{equation}
\!\!\!\left[\begin{array}{c}\sigma^1_{t_n}\\\sigma^2_{t_n}\\\vdots\\\sigma^m_{t_n}\end{array}\right]
\label{Step2}
=\left[\begin{array}{c}\frac{\gamma_{Y_{t_{n-1}}\rightarrow Y_{t_{n}}}(1 )}{\overline\gamma_{Y_{t_{n-1}}\rightarrow Y_{t_{n}}}}\sigma^1_{t_{n}-}\\\frac{\gamma_{Y_{t_{n-1}}\rightarrow Y_{t_{n}}}(2 )}{\overline\gamma_{Y_{t_{n-1}}\rightarrow Y_{t_{n}}}} \sigma^2_{t_{n}-}\\\vdots\\\frac{\gamma_{Y_{t_{n-1}}\rightarrow Y_{t_{n}}}(m )}{\overline\gamma_{Y_{t_{n-1}}\rightarrow Y_{t_{n}}}}\sigma^m_{t_{n}-}\end{array}\right],
\end{equation}
starting at $\!\left[\begin{array}{c}\sigma^1_{t_0}\\\sigma^2_{t_0}\\\vdots\\\sigma^m_{t_0}\end{array}\right]=\left[\begin{array}{c}P(X_0=1)\\P(X_0=2)\\\vdots\\P(X_0=m)\end{array}\right]$.
The unnormalized filter is then
$
\displaystyle \sigma_t(\cdot)=\sum\limits_{i=1}^m \sigma^i_t \delta_i(\cdot)
$.

\section{Conclusions and Highlights}

\subsection{Continuous-time Markov Chains}

Markov chains have many real-world applications such as cruise control systems, queues of customers, financial models, disease spread and population dynamics.
They have also become essential tools in artificial intelligence and in approximating/implementing more complex Markov models.
Their aymptotic mathematical theory is highly developed.
However, general finite-time theory about changing one chain into
another is not fully developed.
We developed a measure change formula and a rejection algorithm
that can be used to change continuous-time Markov chains, add time
dependence or even add hidden states.
In particular, importance sampling, rejection sampling, and Monte Carlo
simulation methods were developed.
These can also be used in the future for parameter estimation and model
learning.

\subsection{Continuous-time Hidden Markov Chains}

Continuous-time Hidden Markov Chains (CTHMM) have become important in
applications like network performance evaluation and 
disease progression tracking.
Many learning approaches have been successfully employed in
the literature, often
based upon the expectation-maximization algorithm.
However, there is no known method to test between local maxima
or between competing hidden models (of possibly different dimensions).
We have provided a Bayes' factor approach, both through a 
particle filter and through a direct method, for solving such
problems through Bayesian model selection.
Further, we have developed new filtering equations for both
the unnormalized and normalized filter as well as direct 
and particle methods to solve them.
These can provide excellent ways to track and predict
disease progression for example.

\subsection{Continuous-time Markov Observation Models}

CTHMMs have the shortcoming that each observation is conditionally
independent of the others. 
This limits their applicability to finance, climate, network security
etc.
To overcome this difficulty, we introduced a new model, the
Continuous-time Markov Observation Models (CMOM).
We proved all our results in this more general context.
Moroever, all applications and solutions discussed apply to CMOM.
In particular, the direct approach to solving the DMZ-type filtering 
equation can be an excellent approach for CTHMM or CMOM when
the transition function of the hidden state is known or can be
well approximated over short time.

\subsection{Mathematics}

Our martingale problem and jump stochastic calculus approach to 
handling Markov chains is unusual.
Our explicit measure-change formula for Markov chains is new
and explicit measure-change formulae are rare. 
There is potential application to simulating (sampling) stochastic
differential equations by Markov chain approximation and simulating (sampling)
the Markov chain approximation with the measure-change formula.
(This would involve approximation but would be far more general
than the Cameron-Martin-Girsanov theorem.)
There are very few Fujisaki-Kallianpur-Kunita or Duncan-Mortensen-Zakai
equations for any type of problem.
Our introduction of the transition count processes into the measure
change and filtering equations is believed to be new.
\begin{funding}
The author gratefully acknowledges support from an NSERC Discovery Grant.
%
\end{funding}

\bibliographystyle{apalike}

\bibliography{PartObsMartProb}

\end{document}